\documentclass{amsart}

\usepackage{amscd,amssymb,amsmath,mathrsfs}
\usepackage[all]{xy}

\newcommand\mylabel[1]{\label{#1}}

\newtheorem{theorem}{Theorem}[section]
\newtheorem*{maintheorem}{Theorem}
\newtheorem{lemma}[theorem]{Lemma}
\newtheorem{proposition}[theorem]{Proposition}
\newtheorem{corollary}[theorem]{Corollary}

\theoremstyle{definition}
\newtheorem{definition}[theorem]{Definition}

\newtheorem*{acknowledgement}{Acknowledgement}

\theoremstyle{remark}

\DeclareFontFamily{U}{wncy}{}
\DeclareFontShape{U}{wncy}{m}{n}{<->wncyr10}{}
\DeclareSymbolFont{mcy}{U}{wncy}{m}{n}
\DeclareMathSymbol{\Sh}{\mathord}{mcy}{"58}
\newcommand{\ZZ}	{\mathbb{Z}}
\newcommand{\QQ}	{\mathbb{Q}}
\newcommand{\RR}	{\mathbb{R}}
\newcommand{\CC}	{\mathbb{C}}

\newcommand{\idealp}	{\mathfrak{p}}

\newcommand{\ideala}	{\mathfrak{a}}

\newcommand{\idealm}	{\mathfrak{m}}

\newcommand{\shA}	{\mathscr{A}}

\newcommand{\shB}	{\mathscr{B}}
\newcommand{\shC}	{\mathscr{C}}

\newcommand{\shE}	{\mathscr{E}}
\newcommand{\shF}	{\mathscr{F}}

\newcommand{\shI}	{\mathscr{I}}

\newcommand{\shL}	{\mathscr{L}}

\newcommand{\calP}	{\mathcal{P}}

\newcommand{\alg}	{{\operatorname{alg}}}

\newcommand{\cH}	{\check{H}}

\newcommand{\et}	{{\text{\rm et}}}
\newcommand{\Et}	{{\text{\rm Et}}}

\newcommand{\Frac}	{\operatorname{Frac}}

\newcommand{\Gal}	{\operatorname{Gal}}

\newcommand{\Hom}	{\operatorname{Hom}}

\newcommand{\id}	{{\operatorname{id}}}

\newcommand{\Ind}	{\operatorname{Ind}}
\newcommand{\ind}	{\operatorname{ind}}

\newcommand{\dirlim}	{\varinjlim}
\newcommand{\invlim}	{\varprojlim}

\newcommand{\lra}	{\longrightarrow}

\newcommand{\Max}	{\operatorname{Max}}
\newcommand{\maxid}	{\mathfrak{m}}
\newcommand{\Min}	{\operatorname{Min}}

\newcommand{\Nis}	{\text{\rm Nis}}

\newcommand{\primid}	{\mathfrak{p}}
\renewcommand{\O}	{\mathscr{O}}

\newcommand{\pr}	{\operatorname{pr}}

\newcommand{\quadand}	{\quad\text{and}\quad}

\newcommand{\ra}	{\rightarrow}

\newcommand{\red}	{{\operatorname{red}}}

\newcommand{\Sing}	{\operatorname{Sing}}

\newcommand{\Spec}	{\operatorname{Spec}}

\newcommand{\TSC}	{\operatorname{TSC}}
\newcommand{\Tor}	{\operatorname{Tor}}

\newcommand{\Zar}	{{\text{\rm Zar}}}

\def\mydate{\number\day\space\ifcase\month \or January\or February\or March\or 
April\or May\or June\or July\or
August\or September\or October\or November\or December\fi \space\number\year}

\DeclareFontFamily{U}{wncy}{}
\DeclareFontShape{U}{wncy}{m}{n}{<->wncyr10}{}
\DeclareSymbolFont{mcy}{U}{wncy}{m}{n}
\DeclareMathSymbol{\Sh}{\mathord}{mcy}{"58}
\begin{document}

\title[Totally separably closed schemes]{Geometry on totally separably closed schemes}

\author[Stefan Schr\"oer]{Stefan Schr\"oer}
\address{Mathematisches Institut, Heinrich-Heine-Universit\"at, 40204 D\"usseldorf, Germany}
%\curraddr{}
\email{schroeer@math.uni-duesseldorf.de}


\dedicatory{Final Version, 21 November 2016}

\begin{abstract}
We prove,   for quasicompact   separated schemes over ground fields, that
\v{C}ech cohomology coincides with sheaf cohomology with respect to  the Nisnevich topology.
This is a partial generalization of Artin's result that for noetherian schemes
such an equality holds with respect to  the \'etale topology, which holds under the assumption that every
finite subset admits an affine open neighborhood (AF-property).
Our key result is that on the absolute integral closure of   separated
algebraic   schemes, the intersection of any two irreducible closed subsets
remains irreducible. We prove this by establishing   general modification and contraction results
adapted to  inverse limits of schemes.
Along the way, we characterize schemes that are acyclic with respect to
various Grothendieck topologies, 
study schemes all local rings of which are strictly henselian, and analyze 
fiber products of strict localizations.
\end{abstract}

\maketitle
\tableofcontents

%===========================================================
\section*{Introduction}
An integral scheme $X$ that is normal and whose function field $k(X)$ is
algebraically closed is called \emph{totally algebraically closed}, or
\emph{absolutely algebraically closed}. 
These notion were introduced by Enochs \cite{Enochs 1968}
 and Artin \cite{Artin 1971} and further studied, for example, in \cite{Hochster 1970a}, \cite{{Hochster 1970b}} and
\cite{Borho; Weber 1971}. 

Artin used such schemes to prove that \v{C}ech cohomology coincides
with sheaf cohomology for any abelian sheaf in the \'etale topos
of a noetherian scheme $X$ that satisfies the \emph{AF-property},
that is, every finite subset admits an affine open neighborhood.
This is a rather large and very useful class of schemes, which includes all schemes that
are quasiprojective of some affine scheme. Note, however, that there are
smooth threefolds lacking this property (Hironaka's example, compare
[31], Appendix B, Example 3.4.1).

In recent years, absolutely closed domains
were studied in connection with  tight closure theory.
Hochster and Huneke \cite{Hochster; Huneke 1992} showed that the absolute
algebraic closure is a big Cohen--Macaulay module in characteristic $p>0$.
This was further extended by Huneke and Lyubeznik  \cite{Huneke; Lyubeznik 2007}.
Schoutens \cite{Schoutens 2004} and Aberbach \cite{Aberbach 2005} used the absolute closure to characterize regularity
for local rings. Huneke \cite{Huneke 2011} also gave a nice survey on absolute algebraic closure.

Other applications include Gabber's rigidity property for   abelian torsion 
sheaves for affine henselian pairs \cite{Gabber 1994}, or Rydh's study of descent questions
\cite{Rydh 2010},
or the proof by   Bhatt and de Jong \cite{Bhatt; de Jong 2014} of the Lefschetz Theorem for local Picard groups.
Closely related ideas occur in the pro-\'etale site, introduced by Bhatt and Scholze \cite{Bhatt; Scholze 2013},
or the construction of   universal coverings for schemes by Vakil and Wickelgren \cite{Vakil and Wickelgren 2011}.
Absolute integral closure in characteristic $p>0$   provides examples of flat
ring extensions that are not a filtered colimit of   finitely presented flat ring extensions,
as Bhatt  \cite{Bhatt 2014} noted. I have used absolute closure to study 
points in the fppf topos \cite{Schroeer 2014}.

For aesthetic reasons, and also to stress the relation to strict localization and   \'etale topology,
I prefer to work with    separable closure instead of   algebraic closure:
If $X_0$ is an integral scheme, 
the \emph{total separable closure}  $X=\TSC(X_0)$ is the integral closure of $X_0$  in some chosen
separable closure of the function field.  
The goal of this paper is to make a systematic study of geometric properties of $X=\TSC(X_0)$,
and to apply it to cohomological questions.
One of our main results is the following rather  counter-intuitive property:

\begin{maintheorem}[compare Thm.\ \ref{no cyclic system}]
Let $X_0$ be separated and of finite type over a ground field $k$, and $X=\TSC(X_0)$.
Then for every pair of closed irreducible subsets $A,B\subset X$, the intersection 
$A\cap B$ remains irreducible.
\end{maintheorem}

Note that for algebraic surfaces, this means that for all closed irreducible subsets $A\neq B$ in $X$,
the intersection $A\cap B$ contains at most one point.
Actually, I conjecture that our result holds true for arbitrary integral   schemes   that are quasicompact 
and separated.

Such geometric properties were crucial for for Artin to establish the equality 
$\cH^p(X_\et,F)= H^p(X_\et,F)$.
He achieved this by assuming  the AF-property. This condition, however, appears to be somewhat  alien to the problem.
Indeed, for locally factorial  schemes, the AF-property is equivalent to quasiprojectivity, according
to Kleiman's proof of the Chevalley conjecture (\cite{Kleiman 1966}, Theorem 3 on page 327).
Recently, this was extended to arbitrary normal schemes by Benoist \cite{Benoist 2013}.
Note, however, that nonnormal schemes may have the AF-property without being quasiprojective,
as examples of Horrocks   \cite{Horrocks 1971} and Ferrand \cite{Ferrand 2003} show.

For schemes $X$ lacking the  AF-property, Artin's argument break down
at two  essential steps: First, he uses the AF-property   to reduce the analysis of fiber products
$\Spec(\O_{X,x}^s)\times_X \Spec(\O^s_{X,y})$ of strict henselizations 
to the more accessible case of integral affine schemes $X=\Spec(R)$, in fact to the situation $A=\ZZ[T_1,\ldots,T_n]$,
which then leads to the join construction  $B=[A_p^h,A_q^h]$ inside some algebraic closure
of the field of fractions $\Frac(A)$ (see \cite{Artin 1971}, Theorem 2.2).
Second, he needs the affine situation to form meaningful semilocal intersection rings $R=A_p^h\cap A_q^h$, 
in order to prove that $B$ has separably closed residue field
(see \cite{Artin 1971}, Theorem 2.5).  

My motivation to study totally separably closed schemes was to bypass the AF-property in Artin's arguments.
In some sense, I was able to generalize half of his reasonings to arbitrary schemes $X$, 
namely those steps that pertain to henselization $\O_{X,x}^h$ for points $x\in X$ rather then strict localizations
$\O_{X,a}^s$ for geometric points $a:\Spec(\Omega)\ra X$, the latter having in addition separably closed residue fields.

In terms of Grothendieck topologies, we thus get results on the \emph{Nisnevich topology}
rather then the \'etale topology. This is one of the more recent Grothendieck topologies, which
was considered in connection with motivic questions.
Recall that the covering families of the  Nisnevich topology on $(\Et/X)$ are those
$(U_\lambda\ra U)_\lambda$ so that each $U_\lambda\ra U$ is \emph{completely decomposed},
that is, over each point  lies at least one point with the same residue field,
and that $\bigcup U_\lambda\ra U$ is surjective. 
We refer to Nisnevich \cite{Nisnevich 1989} and the stacks project \cite{stacks project} for more details.
Note that the Nisnevich topology plays an important role of Voevodsky's theory of sheaves with transfer
and motivic cohomology,
compare \cite{Voevodsky; Suslin; Friedlander 2000}, Chapter 3 and \cite{Mazza; Voevodsky; Weibel 2006}.
Our second main result is:

\begin{maintheorem}[compare Thm.\ \ref{cech = sheaf cohomology}]
Let $X$ be a quasicompact and  separated scheme over a ground field $k$.
Then $\cH^p(X_\Nis,F)=H^p(X_\Nis,F)$ for all abelian Nisnevich sheaves on $X$.
\end{maintheorem}

It is quite sad that my methods apparently need a ground field, in order to use the
geometry of contractions, which lose some of their force over more general ground rings.
Again I conjecture that the result holds true for quasicompact and separated schemes,
even for  the \'etale topology. Indeed, this paper contains several general \emph{reduction steps},
which reveal that it suffices to prove this conjecture merely  for separated  integral $\ZZ$-schemes of finite type.
It seems likely that it suffices that the diagonal
$\Delta:X\ra X\times X$ is affine, rather than closed.

Along the way, it is crucial to characterize schemes that are acyclic with respect to various
Grothendieck topologies.  With analogous results for the Zariski and the \'etale topology, we have:

\begin{maintheorem}[compare Thm.\ \ref{characterization acyclic nisnevich}]
Let $X$ be a quasicompact scheme. Then the following are equivalent:
\begin{enumerate}
\item  $H^p(X_\Nis,F)=0$ for every abelian Nisnevich sheaf $F$ and every $p\geq 1$.
\item Every  completely decomposed \'etale surjection $U\ra X$  admits a section.
\item The scheme $X$ is affine, and each connected component is   local henselian.
\item The scheme $X$ is affine, each irreducible component is local henselian, and the space
$\Max(X)$ is at most zero-dimensional.
\end{enumerate}
\end{maintheorem}

The paper is organized as follows:
In Section \ref{LIMITS}, we review some properties of schemes that are 
stable  by integral surjections. These will be important for many reduction steps that follow.
Section \ref{ESL} contains a discussion of schemes all of whose local rings
are strictly local. Such schemes $X$ have the crucial property that any integral morphism $f:Y\ra X$
with $Y$ irreducible must be injective.
In Section \ref{TOTAL SEPARABLE CLOSURE} we discuss the total separable closure $X=\TSC(X_0)$
of an integral scheme $X_0$, which are the most important examples of schemes that are everywhere
strictly local.

Section \ref{ACYCLIC} contains our characterizations of acyclic schemes with respect to
three Grothendieck topologies, namely Zariski, Nisnevich and \'etale.
In Section \ref{CARTAN--ARTIN}, we introduce technical conditions, namely
the weak/strong Cartan--Artin properties, which roughly speaking means that
the fiber products of strict localizations are acyclic with respect to the Nisnevich/\'etale
topology. Note that such fiber products are almost always non-noetherian.
Section \ref{REDUCTION NORMAL} and \ref{REDUCTION TSC} contain reduction arguments, which basically show that
it suffices to check the Cartan--Artin properties on total separable closures $X=\TSC(X_0)$.
On the latter, it translates into a simple, but rather counter-intuitive geometric condition on
the intersection of irreducible closed subsets.

Section \ref{ALGEBRAIC SURFACES} reveals in the special case of algebraic surfaces this geometric condition indeed holds.
Here one uses that one understands very well which integral curves on a normal surface
are contractible to a point.
The next three sections  prepare the ground to generalize this to higher dimensions:
In Section \ref{FACTS ON SCHEMES}, we collect some facts on noetherian schemes concerning quasiprojectivity
and connectedness of divisors. The latter is an application of Grothendieck's Connectedness Theorem.
In Section \ref{CONTRACTIONS} we show how to make a closed subset contractible on some modification.
In Section \ref{CYCLIC SYSTEMS}, we introduce the technical notion of cyclic systems, which
is better suited to understand contractions in inverse limits like $X=\TSC(X_0)$.
Having this, Section \ref{ALGEBRAIC SCHEMES} contains the Theorem that on total separable closures
there are no cyclic systems, in particular the intersection of irreducible closed subsets
remains irreducible.
The application to Nisnevich cohomology appears in Section \ref{NISNEVICH CECH}.
There are also three appendices, discussing Lazard's observation on connected components of schemes,   H.\ Cartan's 
argument on equality of \v{C}ech and sheaf cohomology, and some results on inductive dimension in general topology
used in this paper.

\begin{acknowledgement}
I wish to thank Johan de Jong for inspiring conversations,
and also the referee for very useful comments, which helped to improve the exposition,
clarify the paper and remove mistakes.
\end{acknowledgement}

%===========================================================
\section{Integral surjections}
\mylabel{LIMITS}
Throughout the paper, integral surjections  and inverse limits   play an important role.
We start  by reexamining these concepts.

Let $X_\lambda$, $\lambda\in L$ be a filtered inverse system of schemes with affine transition morphisms $X_\mu\ra X_\lambda$, 
$\lambda\leq\mu$. Then the corresponding inverse limit  exists
as a scheme. For its construction, one may   tacitly assume that there is a smallest index
$\lambda=0$, and regards the $X_\lambda=\Spec(\shA_\lambda)$ as relatively affine schemes over $X_0$. Then  
$$
X=\invlim(X_\lambda)=\Spec(\shA),\quad \shA=\dirlim\shA_\lambda.
$$
Moreover, the underlying topological space of $X$ is the inverse limit of the underlying topological spaces for $X_\lambda$,
by \cite{EGA IVc},   Proposition 8.2.9. In turn, for each point $x=(x_\lambda)\in X$ the canonical map 
$\dirlim(\O_{X_\lambda,x_\lambda})\ra\O_{X,x}$
is bijective. Consequently, the residue field $\kappa(x)$ is the union of the $\kappa(x_\lambda)$, viewed as subfields.

Recall that a homomorphism $R\ra A$ of rings is   \emph{integral}  if each element in $A$ is the
root of a monic polynomial with coefficients in $R$. 
A morphism of schemes $f:X\ra Y$ is called
integral if there is a affine open covering $Y=\bigcup V_i$ so that $U_i=f^{-1}(V_i)$ is affine
and $A_i=\Gamma(V_i,\O_X)$ is integral as algebra over $R_i=\Gamma(V_i,\O_Y)$.
Note that the underlying topological space of the fibers $f^{-1}(y)$ are profinite.

Of particular interest are those $f:X\ra Y$ that are  integral and surjective.
The following locution will be useful throughout: Let $\calP$ be a class of schemes.
We say that $\calP$ is \emph{stable under images of integral surjections}
if for each integral surjection $f:X\ra Y$ where the domain $X$ belongs to $\calP$,
the range $Y$ belongs to $\calP$ as well.
For example:

\begin{theorem}
\mylabel{permanence affine}
The class  of affine schemes is stable under  images of integral surjections.
\end{theorem}

In this generality, the result is due to Rydh \cite{Rydh 2015}, who established
it even for algebraic spaces.
It generalizes Chevalley's Theorem (\cite{EGA II}, Theorem 6.7.1),
where $Y$  is a noetherian scheme and $f$ is finite   surjective.
See also \cite{Conrad 2007}, Corollary A.2, for the case
that $Y$ is an arbitrary scheme and $f$ is finite surjective.

A scheme $Y$ is called \emph{local} if it  is quasicompact and contains precisely one closed point.
Equivalently  $Y$ is the spectrum of a local ring.
We have the following permanence property:

\begin{proposition}
\mylabel{permanence local}
The class of local schemes is stable under  images of integral surjections
\end{proposition}

\proof
Let $f:X\ra Y$ be an integral surjection, with $X$ local. 
We have to show that $Y$ is local.
According to Theorem \ref{permanence affine}, the scheme $Y$ is affine.
It follows that   $Y$ contains at least one closed point.
Let $y,y'\in Y$ be two closed points. Since $f$ is surjective, there are closed points
$x,x'\in X$ mapping to $y,y'$, respectively. Since $X$ is local, we have $x=x'$, whence $y=y'$.
\qed

\medskip
A scheme $Y$  is called \emph{local henselian} if it is local,
and for every finite morphism $g:Y'\ra Y$, the domain $Y'$ is a sum of local schemes.
Of course, there are only finitely many such  summands, because
$g^{-1}(b)$, where $b\in Y$ is the closed point, contains only finitely many points and  $g$ is a closed map.

\begin{proposition}
\mylabel{permanence local henselian}
The class of   of local henselian scheme 
is stable under images of  integral surjections.
\end{proposition}

\proof
Let  $f:X\ra Y$ be an integral surjection, with $X$ local henselian.
By Proposition \ref{permanence local}, the scheme $Y$ is local. Now let $Y'\ra Y$ be a finite morphism,
and consider the induced finite morphism $X'=Y'\times_YX\ra X$. Then $X'=X'_1\amalg\ldots\amalg X'_n$
for some local schemes $X'_i$, $1\leq i\leq n$. Their images $Y'_i\subset Y'$ are closed, because
$f$ and the induced morphism $X'\ra Y'$ are integral. 
Moreover, these images form a closed covering of $Y$, since $f$ and the induced morphism $X'\ra Y'$ is surjective.
Regarding the $Y_i$ as reduced schemes, the canonical morphism $X'_i\ra Y'_i$ is integral and surjective.
Using Proposition \ref{permanence local} again, we conclude that the $Y'_i$ are local.
It follows that there are at most $n$ closed points on $Y'$. Denote them by $b_1,\ldots,b_m$,
for some $m\leq n$. For each closed point $b_j$, let $C_j\subset Y'$ be the union of those $Y'_i$
that contain $b_j$. Then the $C_j$ are closed  connected subsets,   pairwise disjoint and finite in number. We conclude
that the $C_j$ are the connected components of $Y'$. By construction, each $C_j$ is quasicompact and contains only one
closed point, namely $b_j$, such that $C_j$ is local. 
It follows that $Y'$ is a sum of local schemes, thus $Y$ is local henselian.
\qed

\medskip
A scheme $Y$ is called \emph{strictly local}, if it is local henselian, and the residue field
of the closed point is separably closed. The class of such schemes is not stable under images
of surjective integral morphisms, for example $\Spec(\CC)\ra\Spec(\RR)$.

We say that a class $\calP$ of schemes is stable under  \emph{images of integral surjections  with 
radical residue field extensions} 
if for   every morphism $f:X\ra Y$  that is integral, surjective and whose
field extensions $\kappa(y)\subset\kappa(x)$, $y=f(x)$, $x\in X$ are radical (that is, algebraic and purely inseparable),
and with domain $X$ belonging to  $\calP$, the domain $Y$ also belongs to $\calP$. 

\begin{proposition}
\mylabel{permanence strictly local}
The class of strictly local   schemes is stable  under images of integral surjections   with radical
residue field extensions.
\end{proposition}

\proof
This follows from Proposition \ref{permanence local henselian}, together with the  fact that a field $K$ is separably closed
if it admits a radical extension $K\subset E$ that is separably closed.
\qed

\medskip
Recall that a ring $R$ is called a \emph{pm ring}, if
every prime ideal $\primid\subset R$ is contained in exactly one maximal ideal $\maxid\subset R$,
compare \cite{De Marco; Orsatti 1971}.
Clearly, it suffices to check this for minimal prime ideals, which correspond to the
irreducible components of $\Spec(R)$.
Therefore, we call a scheme $X$ a \emph{pm scheme} if it is quasicompact, and each irreducible component
is local.  These notions will show up in Section \ref{ACYCLIC}. For later use, we observe the following fact:

\begin{proposition}
\mylabel{permanence Gelfand}
The class of pm schemes is stable under images of integral surjections.
\end{proposition}

\proof
Let $f:X\ra Y$ be integral surjective, with $X$ pm. Clearly, $Y$ is quasicompact.
Let $Y'\subset Y$ be an irreducible component. We have to check that $Y'$ is local.
Since $f$ is surjective and closed, there is an irreducible closed subset $X'\subset Y'$ with $f(X')=Y'$.
Replacing $X,Y$ be these subschemes, we may assume that $X,Y$ are irreducible. 
Let $y,y'\in Y$ be two closed points. Choose closed points $x,x'\in X$ mapping to them.
Then $x=x'$, because $X$ is pm, whence $y=y'$.
\qed

%=============================================
\section{Schemes that are everywhere strictly local}
\mylabel{ESL}

We now  introduce a     class of schemes that, in my opinion,  seems   rather natural with respect to the
\'etale topology.
Recall that a \emph{strictly local} ring $R$ is a local ring that is henselian and
has separably closed residue field $\kappa=R/\maxid$. 

\begin{definition}
\mylabel{definition esc}
A scheme $X$ is called \emph{everywhere strictly local} if
the local rings $\O_{X,x}$ are  strictly local, for all $x\in X$.
\end{definition}

Similarly, we call a ring $R$ everywhere strictly local if the
$R_\idealp$ are strictly  local for all prime ideals $\idealp\subset R$.
Note that   strictly local rings are usually not everywhere strictly local.
The relation between these classes of rings 
appears to be similar in nature to the relation between valuation rings and Pr\"ufer   rings
(compare for example \cite{Endler 1972}, \S 11). 
We have the following permanence property:

\begin{proposition}
\mylabel{permanence esl}
If $X$ is everywhere strictly local, the same holds for each quasifinite $X$-scheme.
\end{proposition}

\proof
Let $U\ra X$ be quasifinite, and $u\in U$, with image point $x\in X$. Obviously, the residue field $\kappa(u)$ is  separably closed.
To check that $\O_{U,u}$ is   henselian, we may replace $X$ by the spectrum of $\O_{X,x}$.
According to \cite{EGA IVd}, Theorem 18.12.1 there is an \'etale morphism $X'\ra X$,
a point $u'\in U'=U\times_XX'$ lying over $u\in U$, and an open neighborhood $u'\in V'\subset U'$ so that
the projection $V'\ra X'$ is finite. 
Let $x'\in X'$ be the image of $u'$. After replacing $X'$ by some affine open neighborhood of $x'$,
we may assume that $X'\ra X$ is separated.
There is a section $s:X\ra X'$ with $s(x)=x'$ for the projection $X'\ra X$, because $\O_{X,x}$ is
strictly local, by \cite{EGA IVd}, Proposition 18.8.1. Its image is an open connected component,
according to \cite{EGA IVd}, Corollary 17.9.4, so shrinking further we may assume that $X'=X$.
By \cite{EGA IVd}, Proposition 18.6.8, the local ring $\O_{U,u}=\O_{V',u'}$ is henselian.
\qed

\medskip
Recall that a morphism $f:Y\ra X$ is referred to as  \emph{radical} if it
is universally injective. Equivalently, the map is injective, and
the induced residue field extension $\kappa(x)\subset\kappa(y)$,  are
purely inseparable (\cite{EGA I}, Section 3.7), for all $y\in Y$, $x=f(y)$.
The following geometric property will play a crucial role later: 

\begin{lemma}
\mylabel{integral => radical}
Let $X$ be an everywhere strictly local scheme and $Y$ be an irreducible scheme.
Then every integral morphism $f:Y\ra X$ is radical,   in particular an injective map.
\end{lemma}

\proof
Since the residue fields of $X$ are separably closed and $f$ is integral, it
suffices to check that $f$ is injective.  
Let $y,y'\in Y$ with same image $x=f(y)=f(y')$. Our task is to show $y=y'$.
Replace $X$ by the spectrum of the local ring $R=\O_{X,x}$ and $Y$ by the corresponding fiber product.
Then $Y=\Spec(A)$ becomes affine,
and the morphism of schemes $f:Y\ra X$ corresponds to a homomorphism
of rings $ R\ra A$.
Write $A=\bigcup A_\iota$ as the
filtered union of finite $R$-subalgebras.
Since $R$ is henselian and $A_\iota$ are integral domains, the $A_\iota$ are   local.
Whence $A$ is   local, too. It follows that the closed points $y,y'\in\Spec(A)$ coincide
\qed

\begin{proposition}
\mylabel{local isomorphism}
Let $X$ be everywhere strictly local.
Then every \'etale morphism $f:U\ra X$ is a local isomorphism with respect to the Zariski topology.
\end{proposition}

\proof
Fix a point $u\in U$. We must find an open neighborhood on which $f$ is
an open embedding. Set $x=f(u)$. Consider first the special case that
$f$ admits a section $s:X\ra U$ through $u$.
Since $U\ra X$ is unramified, such a section must be an open embedding by
\cite{EGA IVd}, Corollary 17.4.2. Replacing $U$ by the image of this section,
we reduce to the situation that   $f$ admits a right inverse $s$ that is an
isomorphism. Multiplying $f\circ s=\id_X$ with $s^{-1}$ from the right yields $f=s^{-1}$, which
is an isomorphism.

We now come to the general case. 
Since $\O_{X,x}$ is strictly local, the morphism $U\otimes_X\Spec(\O_{X,x})\ra\Spec(\O_{X,x})$
admits a section through $u\in U$, see \cite{EGA IVd}, Proposition 18.5.11.
According to \cite{EGA IVc}, Theorem 8.8.2, such a section comes from
a  section defined on some open neighborhood of $x\in X$, and the assertion follows.
\qed

\medskip
Given a scheme $X$, we denote by $(\Et/X)$ the   site
whose objects are the \'etale morphisms $U\ra X$, and $(\Zar/X)$ the
site whose objects are the open subschemes $U\subset X$.
In both cases, the covering families $(U_\alpha\ra X)_\alpha$ are  those where the
map $\amalg U_\alpha\ra X$ are surjective.  
In turn, we denote by $X_\et$ and $X_\Zar$ the corresponding
topoi of sheaves within a fixed universe. 
The inclusion functor $i:(\Zar/X)\ra (\Et/X)$ is cocontinuous, which means that
the adjoint $i^*$ on presheaves of the restriction functor $i_*$ preserves the sheaf property.
We thus obtain a morphism of topoi $i:X_\et\ra X_\Zar$.
The following is a direct consequence of the Comparison Lemma (\cite{SGA 4a}, Expos\'e III, Theorem 4.1):

\begin{corollary}
\mylabel{etale = zariski}
Let $X$ be everywhere strictly local. Then the canonical morphism $i:X_\et\ra X_\Zar$
of topoi is an equivalence. In particular, sheaves and cohomology groups for the \'etale
site of $X$ is essentially the same as for the Zariski site.
\end{corollary}

\medskip
We are mainly interested in irreducible schemes. 
Recall that a scheme $X$ is called \emph{unibranch} if it is irreducible,
and the normalization map $X'_\red\ra X_\red$ is bijective (compare \cite{EGA IVa}, Section 23.2.1).
Equivalently, the henselian local schemes $\Spec(\O_{X,x}^h)$ are irreducible, for all $x\in X$,
by \cite{EGA IVd}, Corollary 18.8.16.
We have the following characterization, which
is close to Artin's original definition \cite{Artin 1971}.

\begin{proposition}
\mylabel{characterization tsc}
Let $X$ be irreducible. Then $X$ is everywhere strictly local if and only
if it is unibranch and its function field $k(X)=\kappa(\eta)$ is separably closed.
\end{proposition}

\proof
In light of \cite{EGA IVd} Corollary 18.6.13, the condition is necessary.
To see that it is sufficient, we may assume that $X=\Spec(R)$ is a local integral scheme
and have to check that $R$ is   henselian with separably closed residue field $k=R/\maxid_R$.

Let us start with the latter. Seeking a contradiction,  we assume that $k$ is not separably closed.
Then there is a finite separable field extension $k\subset L$ of degree $d\geq 2$.
By the Primitive Element Theorem, we have $L=k[T]/(f)$ for some
irreducible separable monic polynomial $f\in k[T]$.
Choose a monic polynomial $F(T)$ with coefficients in $R$ reducing to $f(T)$. 
Then $\partial F/\partial T\in R$ is a unit, and it follows
from \cite{Milne 1980}, Chapter I, Corollary 3.6  that   the finite $R$-algebra $A=R[T]/(F)$ is \'etale.
Since the field of fractions $\Omega=\Frac(R)$ is separably closed, we have a decomposition $A\otimes_R\Omega=\prod_{i=1}^d\Omega$.
Set $S=\Spec(A)$, and fix one generic point $\eta_0\in S$. Note that its residue field must be $\Omega$.
Its closure $S_0\subset S$ is a local scheme,  finite over $R$, and thus with residue field  $L$.

Now choose a separable closure $k\subset k'$,    consider the resulting
strict henselization $R'=R^s$, and write $S_0'=S_0\otimes_RR'$ and $S'=S\otimes_RR'$ for
the ensuing faithfully flat base-change.  
According to \cite{SGA 1}, Expos\'e VIII, Theorem 4.1, the preimage  $S_0'\subset S'$ is the closure of the point $\eta_0'\in S'$
lying over $\eta_0\in S$, and thus
$S_0'$ is irreducible, in particular connected.
The closed fiber $S_0'\otimes_{R'}k'=\Spec(L)\otimes_kk'$ is a disjoint union of
$d\geq 2$ closed points, which lie in the same connected component inside $S'$.
But since $R'$ is   henselian, \cite{EGA IVd}, Theorem 18.5.11 (a) ensures
that no two points in the closed fiber $S'\otimes_{R'}k'$ lie in the same
connected component of $S'$, contradiction.
Summing up, the residue field $k=R/\maxid_R$ is separably closed.

It remains to see that $R$ is henselian. 
Let $F\in R[T]$ be a monic polynomial. In light of \cite{EGA IVd}, Theorem 18.5.11 (a')
it suffices to see that the finite flat $R$-algebra $A=R[T]/(F)$ is a product of local rings.
Let $d=\deg(F)=\deg(A)$ be its degree, and consider the $R$-scheme $Y=\Spec(A)$.
Let $Y_1,\ldots,Y_r\subset Y$ be the closures of the connected components of the generic fiber $Y\otimes\Omega$.
Using that $R$ is unibranch, we infer that the closed fibers $Y_i\otimes k$ are local.
By the Going-Down-Theorem for the integral ring extension $R\subset A$, every point
$y\in Y\otimes k$ lies in some $Y_i$. It follows that $Y=\bigcup Y_i$. For each $y\in Y_k$, let $C_y\subset Y$
by the union of all those $Y_i$ containing $y$. Clearly, the $C_y$ are connected and local.
Since the closed fibers of $Y_i$ are local, the $C_y$, $y\in Y_k$ are pairwise disjoint.
In turn, the $C_y$ are the connected components of $Y=\Spec(A)$. It follows that $A$ is a product of local rings.
\qed

%===========================================================
\section{Total separable closure}
\mylabel{TOTAL SEPARABLE CLOSURE}

Let $X$ be an integral scheme. We say that $X$ is \emph{totally separably closed} if
it it is normal, and the function field $k(X)$ is separably closed.
According to Proposition \ref{characterization tsc}, these are precisely the integral normal schemes
that are everywhere strictly local. From this we deduce:

\begin{proposition}
\mylabel{permanence tsc}
If $X$ is totally separably closed, so is every normal closed subscheme $X'\subset X$.
\end{proposition}

Now let $X_0$ be an integral scheme, and choose a  separable closure $F^s$ of the function field $F=k(X_0)$.
We define the \emph{total separable closure} 
$$
X=\TSC(X_0)
$$ 
to be the integral closure of $X_0$ inside $F^s$.
We may regard it as filtered inverse limit: Let $F\subset F_\lambda\subset F^\alg$, $\lambda\in L$ be
the set of intermediate fields that are finite over $F=k(X_0)$, ordered by the inclusion relation,
and let $X_\lambda\ra X_0$
be the corresponding integral closures.  These form a filtered inverse system of schemes
with finite surjective transition maps, and we get a canonical identification
$$
X\lra\invlim X_\lambda,
$$
We tacitly assume that the smallest element in the index set $L$ is denoted by $\lambda=0$.
Note that the   fibers of the map $X\ra X_0$, viewed as a topological space,  are \emph{profinite}.
Recall that such spaces are precisely    totally disconnected compacta,
which are also called \emph{Stone spaces}.

Let me introduce the following notation as a general convention for this paper:
Suppose that  $C\subset X$ is a   closed subscheme. Then the  schematic images 
$C_\lambda\subset X_\lambda$ are closed, and the underlying set is just the image set.
These form an filtered inverse system of schemes,
again with affine transition maps by \cite{EGA II}, Proposition 1.6.2. According to \cite{TG 1-4}, Chapter I, \S4, No. 4, Corollary to Prop.\ 9,  
we get a canonical identification
$C=\invlim C_\lambda$ as sets, and it is easy to see that this is an equality of schemes.
The fiber products $X\times_{X_\lambda}C_\lambda\subset X$, $\lambda\in L$ are closed subschemes
and each one  contains $C$ as a closed subscheme, but is usually much larger.
In fact, one has
$C=\bigcap_{\lambda\in L} (X\times_{X_\lambda}C_\lambda)$
as closed subschemes inside $X$.
Now if $C$ is integral, then its function field $k(C)$ is separably closed. This yields:

\begin{proposition}
If the closed subscheme $C\subset X$ is normal, then $C=\TSC(C_\lambda)$ for each index $\lambda\in L$.
\end{proposition}

\medskip
Now suppose that we have a ground field $k$. Let $X_0$ is an integral $k$-scheme 
and $X=\TSC(X_0)$
its total separable closure, as above. The following   observation reduces
the situation to the case that the ground field is separably closed:
Let $k'$ be the relative algebraic closure of $k$ inside $F^s$, where $F=k(X_0)$.
The scheme $X\otimes_kk'$ is not necessarily integral, but the schematic image $X'\subset X\otimes_kk'$
of the canonical morphism $X\ra X\otimes_kk'$ is.
Note that if the $k$-scheme $X_0$ is algebraic, quasiprojective, or proper,
the respective properties hold for the $k'$-scheme $X'$.  Clearly, the morphism $X\ra X'$ is integral and dominant.
Regarding $k(X)$ as an separable closure of $k(X')$, we 
get a canonical morphism $X\ra \TSC(X')$. 

\begin{proposition}
The canonical morphism $X\ra \TSC(X')$ is an isomorphism.
\end{proposition}

\proof
The scheme $X$ is integral, the morphism in question is integral and birational, and the scheme
$\TSC(X')$ is normal, and the result follows.
\qed

\medskip
Finally, suppose that $X_0$ is an integral \emph{algebraic space} rather than a scheme. Then one may 
define its total separable closure $X=\TSC(X_0)=\invlim  X_\lambda$ in the analogous way. But here
nothing interesting happens:
Indeed, then some $X_\lambda$ is a scheme ([42],
Corollary 16.6.2 when $X_0$ is noetherian, and [53], Theorem B for $X_0$ quasicompact and quasiseparated), 
such that  $\TSC(X_0)$ is a scheme.

%===========================================================
\section{Acyclic schemes}
\mylabel{ACYCLIC}

In this section we study    acyclicity for quasicompact schemes with respect to
the Zariski topology, the Nisnevich topology, and the \'etale topology.  We thus take up the question 
in \cite{SGA 4b}, Expose V, Problem 4.14, to study topoi for which every abelian sheaf has trivial higher cohomology. 

First recall that a topological space  is called \emph{at most zero-dimensional} if  
its topology admits a basis consisting of subsets that are open-and-closed.
Note that this  is in the sense of dimension theory in general topology
(confer, for example, \cite{Pears 1975}), rather than dimension theory in commutative algebra and algebraic geometry.

Given a ring $R$, we write $\Max(R)\subset\Spec(R)$ for the subspace of points
corresponding to maximal ideals. Similarly, we write $\Max(X)\subset X$ for the
set of closed points of a scheme $X$, endowed with the subspace topology.

\begin{theorem}
\mylabel{characterization acyclic zariski}
Let $X$ be a quasicompact scheme. Then the following are equivalent:
\begin{enumerate}
\item We have $H^p(X,F)=0$ for every abelian sheaf $F$   and every $p\geq 1$, where   
cohomology is taken with respect to the Zariski topology.
\item Every surjective local isomorphism $U\ra X$  admits a section.
\item The scheme $X$ is affine, and each connected component is local.
\item The scheme $X$ is affine, each irreducible component is local,
and the space of closed points $\Max(X)$ is at most zero-dimensional.
\item The scheme $X$ is affine, and every  element in $R=\Gamma(X,\O_X)$ is the sum of an idempotent and a unit.
\end{enumerate}
\end{theorem}

\proof 
(i)$\Rightarrow$(ii)
Seeking a contradiction, suppose that some surjective local isomorphism $U\ra X$ does not admit
a section. Since $X$ is quasicompact, there are finitely many affine open subsets $U_1,\ldots,U_n\subset U$
so that each $U_i\ra X$ is an open embedding, and $\amalg U_i\ra X$ is surjective.
Replace $U$ by the direct sum $\amalg U_i$.
Let $F$ be the product of the extension-by-zero sheaves $(f_i)_!(\ZZ_{U_i})$, where
$f_i:U_i\ra X$ are the inclusion maps. Then $\Gamma(X,F)=0$ for all $1\leq i\leq n$.
Consider the \v{C}ech complex
$$
\Gamma(X,F)\lra \prod_{i=1}^n\Gamma(U_i,F)\lra\prod_{i,j=1}^n\Gamma(U_i\cap U_j,F)
$$
The constant section $(1_{U_i})$ in the middle is a cocycle, but not a coboundary, and this
holds true for all refinements of the open covering $X=\bigcup U_i$. In turn, we have $\cH^1(X,F)\neq 0$.
On the other hand, the canonical map $\cH^p(X,F)\ra H^p(X,F)$ from \v{C}ech cohomology to sheaf cohomology is injective and actually bijective 
for $p=1$, whence $H^1(X,F)\neq 0$, contradiction.

(ii)$\Rightarrow$(i)
Let $F$ be an abelian sheaf. Every $F$-torsor becomes trivial on some $U\ra X$ as above.
Since the latter has a section, the torsor is already trivial on $X$.
It follows $H^1(X,F)=0$. In turn, the global section functor $F\mapsto\Gamma(X,F)$
is exact, hence $H^p(X,F)=0$ for all $F$ and  $p\geq 1$.

(ii)$\Rightarrow$(iii)
Choose an affine open covering $X=U_1\cup\ldots\cup U_n$. Using that $\amalg U_i\ra X$ has a section,
we infer that $X$ is quasiaffine, in particular separated.
By the previous implication,  $H^1(X,\shF)=0$ for each quasicoherent sheaf.
According to Serre's Criterion \cite{EGA II}, Theorem 5.2.1, the scheme $X$ is affine.
Next, let $C\subset X$ be a connected component, and suppose $C$ is not local. Choose two different closed points $a_1\neq a_2$ in $C$,
and let $U$ be the   sum of $U_1=X\smallsetminus\left\{a_2\right\}$ and
$U_2=X\smallsetminus\left\{a_1\right\}$.
The ensuing local isomorphism $U\ra X$ allows a section $s:X\ra U$.
Then $s(C)\subset U$ is connected, and intersects both $U_1$ and $U_2$. In turn, $s(C)\subset U_i$ for $i=1,2$.
But $U_1$ and $U_2$ have empty intersection when regarded as subsets of  $U$, contradiction.

(iii)$\Rightarrow$(iv) Let $A\subset X$ be an irreducible component, and $C\subset X$ be its connected component.
Then $A$ is local, because the subset $A\subset C$ is closed and nonempty.
To see that $\Max(X)$ is at most zero-dimensional,  write $X=\Spec(R)$. Let $x\in X$ be a closed point, $x\in U\subset X$ an open neighborhood of
the form $U=\Spec(R_f)$, and $C\subset X$ be the connected component of $x$.
Let $S\subset R$ be the multiplicative system of all idempotents $e\in R$ that are units on $C$.
Then $C=\Spec(S^{-1}R)$, and the localization map $R\ra S^{-1}R$ factors over $R_f$.
In turn, there is some $e\in S$ so that the localization map $R\ra R_e$ factors over $R_f$.
In other words, there is an open-and-closed neighborhood $x\in V$ contained in $U$.
It follows that the subspace $\Max(X)\subset X$ is at most zero-dimensional.

(iv)$\Leftrightarrow$(v) This equivalence   is due to Johnstone \cite{Johnstone 1977}, Chapter V, Proposition in 3.9.

(v)$\Rightarrow$(iii)
Write $X=\Spec(R)$, and let $C=\Spec(A)$ be a connected component, $A=R/\ideala$.
Then every nonunit $f\in A$ is of the form $f=1+u$ for some unit $u\in A^\times$.
In turn, the subset $A\smallsetminus A^\times\subset A$ comprises an ideal, because it coincides with the  Jacobson radical
$J=\left\{f\mid\text{$fg-1\in A^\times$ for all $g\in A$}\right\}$.
Thus $A$ is local.

(iii)$\Rightarrow$(ii)
Let $f:U\ra X$ a surjective local isomorphism. To produce a section, we may assume that
$U$ is a finite sum of affine open subschemes of $X$, in particular of finite presentation over $X$.
Let $x\in X$ be a closed point. 
Obviously, there is a section after base-changing to $C=\Spec(\O_{X,x})$.
Since this is a connected component, we may write $C=\bigcap V_\lambda$, where
the $V_\lambda$ are the open-and-closed neighborhoods of $x$,
compare Appendix A. According to \cite{EGA IVc}, Theorem 8.8.2,
a section already exists over some $V_\lambda$.
Using quasicompactness of $X$, we infer that there is an open-and-closed covering $X=X_1\cup\ldots\cup X_n$
so that sections exists over each $X_i$. By induction on $n\geq 1$, one easily infers that
a section exists over $X$.
\qed

\medskip
Let us call a scheme $X$ \emph{acyclic with respect to the Zariski topology}  if 
it is quasicompact and satisfies the equivalent conditions of the Theorem.
Note that some conditions  in Theorem \ref{characterization acyclic zariski}
already occurred in various other contexts:

Rings  satisfying condition (v), that is, every element is a sum of an idempotent and a unit are also
known as \emph{clean rings}. Such rings 
have been extensively   studied in the realm of commutative algebra (we refer to \cite{Nicholson; Zhou 2005} for an overview).

Rings   for which every prime ideal is contained in only one maximal ideal, one of the conditions
occurring in (iv),  are called \emph{pm rings} or \emph{Gelfand rings}  (they were introduced in  \cite{De Marco; Orsatti 1971}).

Also note that in condition (iv) one cannot remove the assumption that $\Max(X)$ is at most zero-dimensional.
In fact, if $K$ is an arbitrary compact topological space, then the ring $R=\shC(K)$ is pm,
such that its spectrum has local irreducible components; on the other hand,  $\Max(X)=K$ (compare \cite{Gillman; Jerison 1960}, Section 4 for the latter,
and   Theorem 2.11 on p.\ 29 for the former).

\medskip
We now turn to the \emph{Nisnevich topology} on the category $(\Et/X)$.
Recall that a morphism $U\ra V$ of schemes is called \emph{completely decomposed} if for each $v\in V$,
there is a point $u\in U$ with $f(u)=v$ and $\kappa(v)=\kappa(u)$.
The covering families $(U_\alpha\ra U)$ for the Nisnevich topology are those families for which
each $U_\alpha\ra U$ is completely decomposed, and   $\bigcup U_\alpha\ra U$ is surjective.
Sheaves on the ensuing site are referred to as \emph{Nisnevich sheaves}, and we write $H^p(X_\Nis,F)$ for their
cohomology groups, compare \cite{Nisnevich 1989}. 
Note that each  point $x\in X$ yields a point in the sense of topos-theory,
and that the corresponding local ring is the henselization $\O_{X,x}^h$.

\begin{theorem}
\mylabel{characterization acyclic nisnevich}
Let $X$ be a quasicompact scheme. Then the following are equivalent:
\begin{enumerate}
\item  $H^p(X_\Nis,F)=0$ for every abelian Nisnevich sheaf $F$ and every $p\geq 1$.
\item Every  completely decomposed \'etale surjection $U\ra X$  admits a section.
\item The scheme $X$ is affine, and each connected component is   local henselian.
\item The scheme $X$ is affine, each irreducible component is local henselian, and the space
$\Max(X)$ is at most zero-dimensional.
\end{enumerate}
\end{theorem}

The arguments are parallel to the ones for Theorem \ref{characterization acyclic zariski},
and left to  the reader.
For the implication (iv)$\Rightarrow$(iii), one needs the following:

\begin{lemma}
\mylabel{henselian irreducible components}
Let $Y$ be a local scheme.
Then $Y$ is henselian if and only if each irreducible component $Y_i\subset Y$, $i\in I$
is henselian.
\end{lemma}

\proof
The condition is necessary, because $X_i\subset X$ are closed subsets.
Conversely, suppose that each $X_i$ is henselian.
Let $f:X\ra Y$ be a finite morphism, and let $a_1,\ldots,a_n\in X$ be the closed points. 
Each of them maps to the closed point $b\in Y$.
The closed subsets $X_i=f^{-1}( Y_i)$ contain  $a_1,\ldots,a_n$ and are finite over $X_i$,
whence the canonical map $\amalg_{j=1}^n\Spec(\O_{X_i,a_j})\ra X_i$ is an isomorphism.
For each $1\leq j\leq n$, consider the subset
$$
C_j=\bigcup_{i\in I} \Spec(\O_{X_i,a_j})\subset X.
$$
Then the $C_j\subset X$ are stable under generization and  contain a single closed point $a_j$,
thus $C_j=\Spec(\O_{X,a_j})$. Furthermore, the $C_j$ form a 
a partition of $X$. It is not a priori clear that $C_j\subset X$ is closed if the index
set $I$ is infinite. But this nevertheless holds: Write $X=\Spec(A)$, where $A$ 
is a semilocal ring. Let $\idealm_j\subset A$ be the maximal ideal corresponding to $a_j\in X$.
Since the $C_j$ are pairwise disjoint, the ideals $\idealm_j$ are pairwise coprime. By the Chinese Reminder Theorem,
the canonical map $A\ra\prod_{j=1}^nA_{\idealm_j}$ is bijective. Whence $C_j\subset X$
are closed, and $X$ is the sum of local schemes.
\qed
 
\medskip
We finally come to the \'etale topology. This is the Grothendieck topology on the category
$(\Et/X)$ whose covering families $(U_\alpha\ra U)$ are those families with
$\bigcup U_\alpha\ra U$ surjective.
Sheaves on this site are called \emph{\'etale sheaves}, and we write $H^p(X_\et,F)$ for their cohomology groups.
In the following, the equivalence (ii)$\Leftrightarrow$(iii) is due to Artin \cite{Artin 1971}:

\begin{theorem}
\mylabel{characterization acyclic etale}
Let $X$ be a quasicompact scheme. Then the following are equivalent:
\begin{enumerate}
\item We have $H^p(X_\et,F)=0$ for every abelian \'etale sheaf $F$ and every $p\geq 1$.
\item Every \'etale surjective  $U\ra X$  admits a section.
\item The scheme $X$ is affine, and each connected component is   strictly local.
\item The scheme $X$ is acyclic, each irreducible component is strictly local, and the space
$\Max(X)$ is at most zero-dimensional.
\end{enumerate}
\end{theorem}

The remaining arguments are parallel to the ones for Theorem \ref{characterization acyclic zariski},
and left to  the reader.

I conjecture that the class of schemes that are acyclic with respect to either the Zariski or the Nisnevich topology
are stable under images of integral surjections. 
Indeed, if  $f:X\ra Y$ is  integral and surjective, with $X$ acyclic,
then $Y$ is affine by Theorem  \ref{permanence affine}, and it follows from Proposition \ref{permanence local} and \ref{permanence local henselian} that
each irreducible component is local or local henselian, respectively. 
It only remains to check that $\Max(Y)$ is at most zero-dimensional, and here lies the problem:
We have a closed continuous surjection $\Max(X)\ra\Max(Y)$,
but   cannot conclude from this that the image  is at most zero-dimensional, in light
of the existence of dimension-raising maps (confer \cite{Pears 1975}, Chapter 7, Theorem 1.8).
The following weaker statement will suffice for our applications:

\begin{proposition}
\mylabel{covering acyclic zariski}
Let $Y$ be a scheme, and $Y=Y_1\cup\ldots \cup Y_n$ be a closed covering.
If each $Y_i$ is acyclic with respect to the Zariski topology, so is $Y$.
\end{proposition}

\proof
As discussed above, it remains to show that the  space of closed points 
$\Max(Y)$ is at most zero-dimensional.
The canonical morphism $Y_1\amalg\ldots\amalg Y_n\ra Y$ is surjective and integral, whence $Y$ is a pm scheme
by Proposition \ref{permanence Gelfand}. Note that this ensures that the topological space $Y$ is \emph{normal}, that is, 
disjoint closed subsets can be separated by disjoint open neighborhoods, 
according to \cite{Carral; Coste 1983}, Proposition 2.
Furthermore, the subspace  $\Max(Y)$ is \emph{compact}, which means quasicompact and Hausdorff 
(\cite{Schwartz 2011}, Corollary 4.7 together with Proposition 4.1).
Since our covering is closed, we have $\Max(Y)=\Max(Y_1)\cup\ldots\cup\Max(Y_n)$, and this is again a closed covering.

To proceed, we now use the notion of 
small and large inductive dimension from general topology.
Since this material is perhaps not so well-known outside general topology,
I have included a brief discussion in Appendix C.
The relevant results are tricky, but rely only on basic notations of topology.
Proposition \ref{ZeroDimensionalCovering}, which is a consequence of a general sum theorem on the large inductive dimension,
ensures that the compact space $\Max(Y)=\Max(Y_1)\cup\ldots\cup\Max(Y_n)$ is at most zero-dimensional.
\qed

%===========================================================
\section{The Cartan--Artin properties}
\mylabel{CARTAN--ARTIN}

Let $X$ be a  topological space and $F$ be an abelian sheaf. Then there is a spectral sequence
$$
E_2^{pq}=\check{H}^p(X,\underline{H}^q(F)) \Longrightarrow H^{p+q}(X,F)
$$
computing sheaf cohomology in terms of \v{C}ech cohomology, where
$\underline{H}^q(F)$ denotes the presheaf $U\mapsto H^p(U,F)$.
A result of H.\ Cartan 
using this spectral sequence tells us that \v{C}ech and sheaf cohomology
coincide on spaces admitting a basis $\shB$ for the topology  that is stable under finite intersections
and satisfies $H^p(U,F)=0$ for all  $U\in\shB$, $p\geq 1$, $F$. 

The same holds, cum grano salis,  for arbitrary sites (compare   the discussion in Appendix B).
For the \'etale site of a scheme $X$, there
seems to be no   candidate for such a basis $\shB$ of open sets. 
However, there is the following substitute:
Let $a=(a_1,\ldots,a_n)$ be a sequence of geometric points $a_i:\Spec(\Omega_i)\ra X$.
Following Artin, we define 
$$
X_a=X_{a_1,\ldots,a_n}=\Spec(\O_{X,a_1}^s)\times_X\ldots\times_X\Spec(\O_{X,a_n}^s).
$$
The schemes $X_a$  can be viewed as inverse limits
of  finite intersections of   small and smaller   neighborhoods, and therefore appears to be
a good substitute for the $U\in\shB$. It was M.\ Artin's insight 
\cite{Artin 1971}   that  for the collapsing of the Cartan--Leray spectral sequence in the \'etale topology
it indeed suffices to verify that  the inverse limits $X_a$ are acyclic.

In order to treat the Nisnevich topology, we use the following:
If $x=(x_1,\ldots,x_n)$ is a sequence of ordinary points $x_i\in X$, we likewise set
$$
X_x=X_{x_1,\ldots,x_n}=\Spec(\O_{X,x_1}^h)\times_X\ldots\times_X\Spec(\O_{X,x_n}^h),
$$
using the henselizations rather then the strict localizations.
Usually, one writes $\overline{x}_i:\Spec(\Omega_i)\ra X$ for   geometric points with image point
$x_i\in X$. In order to keep notation simple, and since we are mainly interested in 
geometric points, I prefer to write $a_i $ for geometric points.
Although this is slightly ambiguous, it should   cause no confusion.
Let us now introduce the following terminology:

\begin{definition}
\mylabel{definition cap}
We say that the  scheme $X$ has the \emph{weak Cartan--Artin property} if
for each sequence $x=(x_1,\ldots,x_n)$  of points on $X$, the scheme
$X_x$ is acyclic with respect to the Nisnevich topology.
 We say that $X$ has the \emph{strong Cartan--Artin property} if for each sequence
$a=(a_1,\ldots,a_n)$ of geometric points on $X$, the scheme
$X_a$ is acyclic with respect to the \'etale topology.
\end{definition}

In the latter case, this means that the schemes $X_a$ are affine, their irreducible components
are strictly local, and the subspace $\Max(X_a)\subset X_a$ of closed points is at most zero-dimensional.
In the former case, however, there is no condition on the residue fields of the closed points.
Our interest in this property comes from the following result, which was proved by Artin \cite{Artin 1971}
for the strong Cartan--Artin property. For the weak Cartan--Artin property, the arguments
are virtually the same.

\begin{theorem}[\textbf{M.\ Artin}]
Let $X$ be a noetherian scheme.
If it has the  strong Cartan--Artin property, then $\check{H}^p(X_\et,\underline{H}^q(F))=0$ for all
abelian \'etale sheaves $F$ and all $p\geq 0$, $q\geq 1$. 
If it has the weak Cartan--Artin property, then $\check{H}^p(X_\Nis,\underline{H}^q(F))=0$ for all
abelian Nisnevich sheaves $F$ and all $p\geq 0$, $q\geq 1$. 
In particular, \v{C}ech cohomology coincides
with sheaf cohomology.
\end{theorem}

To avoid tedious repetitions, we will work in both cases with the schemes $X_a$ constructed
with the strict localizations. This is permissible, in light of the following observation:

\begin{proposition}
\mylabel{characterization wcap}
The scheme $X$ has the weak Cartan--Artin property if and only if for each sequence 
$a=(a_1,\ldots,a_n)$ of geometric points on $X$, the scheme $X_a$ is acyclic with respect to the
Nisnevich topology.
\end{proposition}

\proof
Let $x_i\in X$ be the image of the geometric point $a_i$, and set $x=(x_1,\ldots,x_n)$.
Consider the canonical morphism $X_a\ra X_x$. The rings $\O_{X,a}^s$ are filtered direct limits
of finite \'etale $\O_{X,x}^h$-algebras, whence , the scheme $X_a$ is an filtered inverse limit $X_a=\invlim V_\lambda$
of $X_x$-schemes whose structure morphism $V_\lambda\ra X_x$ are \'etale coverings.
In particular, the morphism $f:X_a\ra X_x$ and the projections $X_a\ra V_\lambda$  are a flat integral   surjections.

We now make the following observation:
Suppose $B\subset X_a$ and $C\subset X_x$ are connected components, with $f(B)\subset C$.
We claim that the induced map $f:B\ra C$ is surjective.
To see this, let $B_\lambda\subset V_\lambda$ be the connected components containing the images of $B$.
These $B_\lambda$ form a filtered inverse system of connected affine schemes with $B=\invlim B_\lambda$,
and the structure maps $B_\lambda\ra X_x$ factor  over $C$. Since $C$ is connected
and the induced projections $V_\lambda\times_{X_x}C\ra C$ are \'etale coverings, the preimages $V_\lambda\times_{X_x}C$
are sums of finitely many connected components, each being an \'etale covering of $C$.
This follows, for example, from the fact that the category of finite Galois coverings of
a connected scheme is a Galois category  admitting a fiber functor. This was established for locally noetherian
schemes in \cite{SGA 1}, Expose V; for the general case see \cite{Lenstra 1985}, Theorem 5.24.
One of the connected components must be $B_\lambda$, and we conclude that the $B_\lambda\ra C$ are surjective.
Likewise, the transition maps $B_\lambda\ra B_\mu$, $\lambda\geq \mu$ must be surjective, and we infer
that $B=\invlim B_\lambda\ra C$ is surjective.

Now suppose that $X_a$ is acyclic with respect to the Nisnevich topology.
Then it is affine by Theorem \ref{characterization acyclic nisnevich}. It follows from Proposition \ref{permanence affine}
that $X_x$ is affine. 
Let $C\subset X_x$ be a connected component. We have to show that $C$ is local.
Since $f:X_a\ra X_x$ is surjective, there is a connected component $B\subset X_a$ with $f(B)\subset C$,
and we saw in the preceding paragraph that $f:B\ra C$ is an integral surjection.
Since $X_a$ is acyclic, the scheme $B$ is local henselian. The same then holds for $C$, by 
Proposition \ref{permanence local henselian}.

Conversely, suppose that $X_x$ is acyclic.
Then $X_x$ is affine, and the same holds for $X_a$, because $f$ is an affine morphism.
Let $B\subset X_a$ be a connected component, and $C\subset X_x$ the connected component
containing $f(B)$. Then $B$ is local, and we saw above that $f:B\ra C$ is an integral surjection.
Moreover, the $B_\lambda$ are local, and the transition maps are local, whence $B=\invlim B_\lambda$ is local.
\qed

%===========================================================
\section{Reduction to  normal schemes}
\mylabel{REDUCTION NORMAL}

The goal of this section is to reduce  checking the Cartan--Artin properties to the case of irreducible, or even normal schemes.
Artin \cite{Artin 1971} bypassed this by assuming the AF-property, together with the fact that any affine scheme is a 
subscheme of some integral scheme. Without assuming the AF-property,  or knowing that any scheme is a subscheme
of some integral scheme, a   different approach is necessary.

Let $X$ be a scheme.
We now take a closer look at the functoriality of  the   $X_a$ with respect to $X$.
Let $f:X\ra Y$ be a morphism of schemes, and $a=(a_1,\ldots,a_n)$ a sequence of geometric points on $X$,
and $b=(b_1,\ldots,b_n)$ the sequence of geometric image points on $Y$.
We then obtain an induced map between strictly local rings, and thus a morphism $X_a\ra Y_b$.

\begin{lemma}
\mylabel{induced separated}
Let $\calP$ be a class of morphisms that is stable under fiber products, and that contains all closed embeddings.
If $f:X\ra Y$ is separated, and the induced morphisms $\Spec(\O^s_{X,a_i})\ra\Spec(\O^s_{Y,b_i})$
belong to $\calP$, then   $X_a\ra Y_b$ belongs to $\calP$.
\end{lemma}

\proof
By assumption, the induced morphism
$$
\prod_{i=1}^n\Spec(\O_{X,a_i}^s)\lra \prod_{i=1}^n\Spec(\O_{Y,b_i}^s)
$$
belongs to $\calP$, where the   products designate  fiber products over $Y$.
It remains to check the following: If $U,V$ are two $X$-schemes, than
the morphism
$U\times_X V\ra U\times_Y V$
is a closed embedding. This indeed holds by \cite{EGA I}, Proposition 5.2, because  $f:X\ra Y$ is separated.
\qed

\begin{proposition}
\mylabel{induced integral}
Notation as above. If $f:X\ra Y$ is integral then  the same   holds for $X_a\ra Y_b$.
If moreover $f:X\ra Y$ is   finite,  a closed embedding, radical or with radical field extensions,   the respective property
holds for $X_a\ra Y_b$.
\end{proposition}

\proof
In light of Lemma \ref{induced separated}, it suffices to treat the case 
that the sequence $a=(a_1,\ldots,a_n)$ consists of a single geometric point,
which by abuse of notation we also  denote by $a$. 
Write $X=\invlim X_\lambda$ as a filtered inverse limit of finite $Y$-schemes $X_\lambda$,
and let $a_\lambda$ be the geometric point on $X_\lambda$ induced by $a$.
Let $x\in X$ and $x_\lambda\in X_\lambda$ be their respective images.
Then $\O_{X,x}=\dirlim\O_{X_\lambda,x_\lambda}$, and it follows from \cite{EGA IVd}, Corollary 18.6.14
that $\O^s_{X,a}=\dirlim\O_{X_\lambda,a_\lambda}^s$.
On the other hand, the base change $X_\lambda\times_Y\Spec(\O_{Y,b}^s)$ is a finite sum of   strictly local
schemes. Let $C_\lambda$ be the connected component corresponding to $a_\lambda'=(a_\lambda,b)$.
It follows from \cite{EGA IVd}, Proposition 18.6.8 that $\O_{X_\lambda,a_\lambda}^s=\O_{C_\lambda,a_\lambda'}$.
In turn, $\O_{X,a}^s$ is a filtered inverse limit of finite $\O_{Y,b}^s$-algebras, thus integral.
From this description, the additional assertions follow as well.
\qed

\medskip
Our   goal now is to reduce checking the  Cartan--Artin properties to the case of integral normal schemes.
Recall that $a=(a_1,\ldots,a_n)$ is a fixed sequence of geometric points on $X$.

\begin{proposition}
\mylabel{irreducible closed covering}
Suppose there are only finitely many irreducible components $X_j\subset X$, $1\leq j\leq r$ that contain
all geometric points $a_1,\ldots,a_n$. Then the closed subschemes
$$
X_{j,a}=(X_j)_a=\prod_{i=1}^n\Spec(\O_{X_j,a_i}^s),\quad 1\leq j\leq r 
$$
form  a closed covering of $X_a$. Here the product means fiber product over $X$.
\end{proposition}

\proof
Let $x_i\in X$ be the  images of the geometric points $a_i$. Clearly, the structure morphism $X_a\ra X$ factors through
$\bigcap_{i=1}^n \Spec(\O_{X,x_i})$,
and this intersection is the set of points $x\in X$ containing all $x_1,\ldots,x_n$ in their closure.
Thus the inclusion 
$\bigcap_{i=1}^n \Spec(\O_{X',x_i})\subset\bigcap_{i=1}^n \Spec(\O_{X,x_i})$
is an equality, where $X'=X_1\cup\ldots\cup X_r$.
It then follows that the closed embedding $X'_a\subset X_a$ is bijective.
We thus may assume that $X=X'$.

Let $\eta\in X_a$ be a generic point. The structure morphism $X_a\ra X$ is flat.
By going-down, the point $\eta$ maps to the generic point $\eta_j\in X_j$ for some $1\leq j\leq r$.
In turn, it is contained in $(X_j)_a=X_a\times_XX_j$. 
Consequently the $(X_j)_a\subset X$ form a covering, 
because they are contain every generic point and are closed subsets.
\qed

\begin{proposition}
\mylabel{irreducible Cartan--Artin}
Suppose the irreducible components of $X$ are locally finite, and
that the diagonal $ X\ra X\times X$ is an affine morphism.
Then $X$ has the weak/strong Cartan--Artin property if and only if the respective property holds for each irreducible component,
regarded as an integral scheme.
\end{proposition}

\proof
Let $a=(a_1,\ldots,a_n)$ be a sequence of geometric points on $X$, and  denote by
$X_1,\ldots,X_r\subset X$   the irreducible components over which all $a_i$ factors.
According to \ref{irreducible closed covering}, we have a closed covering $X_a=(X_1)_a\cup\ldots\cup(X_r)_a$.
Now Proposition \ref{covering acyclic zariski} ensures that the space of closed points $\Max(X_a)$ is at most zero-dimensional.
Let $S=(X_1)_a\amalg\ldots\amalg(X_r)_a$ be the sum. Then $S\ra X_a$ is integral and surjective.
By the assumption on the diagonal morphism of $X$, the scheme
$S$ is affine, and the scheme $X_a$ must be affine as well, by Theorem \ref{permanence affine}.
Finally, let $C\subset X_a$ be an irreducible component. Choose some $1\leq j\leq r$ with $C\subset (X_j)_a$.
Then $C$ is an irreducible component
of $(X_j)_a$.  Since $(X_j)_a$ is acyclic, the scheme $C$ is strictly local, according to Theorems 
\ref{characterization acyclic nisnevich} and \ref{characterization acyclic etale},
and it follows that $X_a$ is acyclic.
\qed

\medskip
Concerning the Cartan--Artin properties, we   now may restrict our attention to integral schemes.
Our next task is to reduce to normal schemes.      
Suppose   that $f:X\ra Y$ is a finite birational morphism between integral schemes,
and let $b=(b_1,\ldots,b_n)$ be a sequence of geometric points in $Y$.  As a shorthand,  we write 
$$
f^{-1}(b)=f^{-1}(b_1)\times\ldots\times f^{-1}(b_n)
$$
for the finite set of sequences $a=(a_1,\ldots,a_n)$ of geometric points with $f(a_i)=b_i$.
Note that one may regard $f^{-1}(b)$ as the product of the points in the finite schemes $X\times_Y\Spec\kappa(b_i)$.
For each $a\in f^{-1}(b)$, we obtain an induced morphism $X_a\ra Y_b$, which is finite.
Its image is denote, for the sake of simplicity, by  $f(X_a)\subset Y_b$, which is a closed
subset.
Now recall that an integral scheme $X$ is called \emph{geometrically unibranch} if the normalization map $X'\ra X$
is radical. This   obviously holds if $X$ is already normal.

\begin{proposition}
\mylabel{normalization closed covering}
Assumptions as above.  Suppose that $X$ is geometrically unibranch.
Then the finite morphisms $X_a\ra Y_b$ are radical,
and the   closed subsets 
$$
f(X_a)\subset Y_b,\quad a\in f^{-1}(b)
$$
form a closed covering of $Y_b$.
\end{proposition}

\proof
First, let us check that $X_a\ra Y_b$ is radical.
In light of Proposition \ref{induced integral}, it suffices to treat the case that the sequence $a=(a_1,\ldots,a_n)$
consists of a single geometric point.
Let $x\in X$, $y\in Y$ be the points corresponding to the geometric points $a=a_1$ and $b=b_1$, respectively.
Set $R=\O_{Y,y}$ and write $X\times_Y\Spec(\O_{Y,y})=\Spec(A)$. Then $A$ is a semilocal, integral,
and a finite $R$-algebra. According to \cite{EGA IVd}, Proposition 18.6.8, we have
$A\otimes_R R^s=A^s$. This is a finite $R^s$-algebra, whence splits into $A^s=B_1\times\ldots\times B_m$,
where the factors are local and correspond to the geometric  points in $f^{-1}(y)$. In fact, the factors
are the strict localizations of $X$ at the these geometric points, and $\O_{X,a}$ is one of them. Since $X$ is geometrically unibranch,
all factors $B_i$ are integral, according to \cite{EGA IVd}, 18.6.12.

As $f:X\ra Y$ is birational, the inclusion $R\subset A$ becomes bijective  
after localizing with respect to the multiplicative system $S=R\smallsetminus 0$.
In turn, the inclusion on the left
$$
S^{-1}R^s\subset S^{-1}A^s=S^{-1}(A\otimes_RR^s)=S^{-1}A\otimes_{S^{-1}R}S^{-1}R^s=S^{-1}R^s
$$
is bijective, such that we have a canonical identification $S^{-1}A^s=S^{-1}R^s$.
Let $R_i\subset B_i$ be the image of the composite map $R^s\ra B_i$, which is also a residue class ring of
$R^s$. Hence the $R_i\simeq R^s$ are integral strictly local rings, in particular geometrically unibranch,
and $R_i\subset B_i$ induces a bijection on fields of fractions.
It follows that $\Spec(B_i)\ra\Spec(R_i)$ are radical. 
It follows that  $\Spec(\O^s_{X,a})\ra\Spec(\O^s_{Y,b})$ is  radical.

It remains to check that $Y_b=\bigcup_af(X_a)$. Let $\eta\in Y_b$ be a generic point.
The images $\eta_i\in\Spec(\O_{Y,b_i})$ are generic points as well, 
because the projections $Y_b\ra\Spec(\O^s_{Y,b_i})$ are flat, whence satisfy going-down.
We saw in the preceding paragraph that the generic points in $\Spec(\O^s_{Y,b_i})$ correspond
to the   points in $f^{-1}(b_i)$.
Let $a_i\in f^{-1}(b_i)$ be the point corresponding to $\eta_i$, and form the sequence $a=(a_1,\ldots,a_n)$.
Using the $f:X\ra Y$ is birational, we infer that $\eta$ lies in the image of $X_a\ra Y_b$. 
\qed

\medskip
If $Y$ is an integral scheme, its normalization map $f:X\ra Y$ is integral, but not necessarily
finite. However, the latter holds for the so-called \emph{japanese schemes}, confer \cite{EGA IVa}, Chapter 0, \S23.

\begin{proposition}
\mylabel{normalization Cartan--Artin}
Suppose that $Y$ is an integral scheme whose normalization morphism $f:X\ra Y$ is finite.
If $X$ has the weak/strong Cartan--Artin property, then the respective property holds for  $Y$.
\end{proposition}

\proof
Let $b=(b_1,\ldots,b_n)$ be a sequence of geometric points on $Y$.
One argues as for the proof of Proposition \ref{irreducible Cartan--Artin}.
 Using that $f:X_a\ra Y_b$ are radical, one sees that
the closed points on $Y_b$ have residue fields that are separably closed.
Whence $Y_b$ is acyclic with respect to the Nisnevich/\'etale topology.
\qed

%===========================================================
\section{Reduction to TSC schemes}
\mylabel{REDUCTION TSC}

I conjecture that the class of schemes having the weak Cartan--Artin property
is stable under images of integral surjections. The goal of this section
is to establish a somewhat weaker variant that suffices for our applications.
Recall that a scheme $X$ is called \emph{geometrically unibranch}
if it is irreducible, and  the normalization morphism $Y\ra X_\red$ is radical.
Equivalently, the strictly local rings $\O_{X,a}^s$ are irreducible, for all
geometric points $a$ on $X$, according to \cite{EGA IVd}, Proposition 18.8.15.

\begin{theorem}
\mylabel{image weak cartan--artin}
Let $f:X\ra Y$ be an   integral surjection between 
schemes that are  geometrically unibranch and have affine diagonal.
If $X$ has the weak Cartan--Artin property, so does  $Y$.
\end{theorem}

This relies on a precise understanding of the connected components inside the $X_a$,
for sequences $a=(a_1,\ldots,a_n)$ of geometric points. We start by collecting some useful facts.

\begin{proposition}
\mylabel{connected components integral}
Suppose $X$ is   geometrically unibranch and quasiseparated. 
Then the connected components of $X_a$ are irreducible.
\end{proposition}

\proof
Clearly, the scheme 
$$
X_a=X_{a_1,\ldots,a_n}=\Spec(\O_{X,a_1}^s)\times_X\ldots\times_X\Spec(\O_{X,a_n}^s)
$$ 
is a filtered inverse limit of \'etale $X$-schemes that are quasicompact and quasiseparated.
Let $U\ra X$ be an \'etale morphism, with $U$ affine.
Since $U\ra X$ satisfies going-down, each generic point on $U$ maps to the generic point in $X$.
Thus $U$ contains only finitely many irreducible components, because $U\ra X$ is quasifinite.
Whence for each connected component $U'\subset U$, there
is a sequence of irreducible components $U_1,\ldots,U_n\subset U$, possibly with repetitions,
so that $U'=\bigcup U_i$ and $U_i\cap U_{i+1}\neq\emptyset$. Choose such a sequence with $n\geq 1$ minimal.
Seeking a contradiction, we suppose  $n\geq 2$. Choose
 $u\in U_1\cap U_2$. Let $x\in X$ be its image. Clearly, $\O_{X,x}\ra\O_{U,u}$ becomes
bijective upon passing to strict henselizations. Since $\O_{X,x}^s$ is irreducible
and $\O_{U,u}\subset\O_{U,u}^s$ satisfies going-down, we conclude that $\O_{U,u}$ is irreducible,
contradiction. The upshot is that the connected components of $U$ are irreducible.
It remains to apply the Lemma below with $X_a=V$.
\qed

\begin{lemma}
\mylabel{limits integral}
Let $V=\invlim V_\lambda$ be a filtered inverse limit with flat transition maps of schemes
$V_\lambda$ that are quasicompact and quasiseparated, and whose connected components  irreducible.
Then each connected component of $V$ is irreducible.
\end{lemma}

\proof
Let $C\subset V$ be a connected component. Then the closure  
of the images $\pr_\lambda(C)\subset V_\lambda$ are connected, whence are contained
in a connected component $C_\lambda\subset V_\lambda$. These connected components
form an inverse system, and we have canonical maps
$C\ra\invlim C_\lambda\subset\invlim V_\lambda=V$.
The inclusion is a closed embedding, and $\invlim C_\lambda$ is connected, since the $C_\lambda$ are
quasicompact and quasiseparated.
In turn $C=\invlim C_\lambda$.  By assumption, the closed subsets $C_\lambda\subset V_\lambda$ are irreducible.
According to a result of Lazard discussed in Appendix A,
the connected component $C_\lambda\subset V_\lambda$ is the intersection
of its neighborhoods that are open-and-closed.
This ensures that the transition maps $C_\mu\ra C_\lambda$ remain flat. If follows that the projections
$C\ra C_\lambda$ are flat as well. By going-down, each generic point $\eta\in C$ maps
to the generic point $\eta_\lambda\in C_\lambda$, whence $\eta=(\eta_\lambda)$ is unique,
and $C$ must be irreducible.
\qed

\medskip
Suppose that $X$ is  geometrically unibranch and quasiseparated,
such that the connected components of $X_a$ coincide with
the irreducible components. We may describe the subspace $\Min(X_a)\subset X_a$ of generic points
as follows: Let $K=k(X)$ be the function field, and $\O_{X,a_i}^s\subset L_i$ be the field
of fractions of the strictly local ring. Since $f:X_a\ra X$ is flat and thus satisfies going-down,
each generic point of $X_a$ maps to the generic point $\eta\in X$. Since $X_a\ra X$ is an inverse limit of quasifinite
$X$-scheme, each point in $f^{-1}(\eta)$ is indeed a generic point of $X_a$. Thus:

\begin{proposition}
\mylabel{tensor product fields}
Assumptions as above. Then the subspace of generic points $\Min(X_a)\subset X_a$ is canonically identified
with $\Spec(L_1\otimes_K\ldots\otimes_K L_n)$.
\end{proposition}

In particular, the space $\Min(X_a)$ is profinite. We may describe it in terms of Galois theory
as follows: Choose a separable closure $K\subset K^s$ and embeddings $L_i\subset K^s$.
Let $G=\Gal(K^s/K)$ be the Galois group, and set $H_i=\Gal(K^s/L_i)$. Then $G$ is a profinite
group, the $H_i\subset G$ are closed subgroups, and the quotients $G/H_i$ are profinite. Then:

\begin{proposition}
\mylabel{quotient product homogeneous}
Assumptions as above. Then the space of generic points $\Min(X_a)$ is homeomorphic
to the orbit space for the canonical left $G$-action on the    space $G/H_1\times\ldots\times G/H_n$.
\end{proposition}

\medskip
\emph{Proof of Theorem \ref{image weak cartan--artin}.}
Suppose that $X,Y$ are  geometrically unibranch, with affine diagonal, and let
$f:X\ra Y$ be a   integral surjection. Assume   $X$ satisfies the
weak Cartan--Artin property. We have to show that the same holds for $Y$. Without restriction,
it suffices to treat the case that $Y,X$ are integral.

Let $b=(b_1,\ldots,b_n)$ be a sequence of geometric points on $Y$, and $C\subset Y_b$ be
a connected component. Our task is to check that $C$ is local henselian.
According to Proposition \ref{connected components integral}, the scheme $C$ is integral. Let $\eta\in C$ be the generic point.
Since $f$ is surjective and integral, we may lift   $b$ to a sequence of geometric points $a$ on $X$.
Since $X,Y$ are geometrically unibranch, the strictly local rings $\O_{X,a_i}^s,\O_{Y,b_i}^s$ remain integral,
and we have a commutative diagram
$$
\begin{CD}
\O_{X,a_i} @>>> K_i\\
@AAA @AAA\\
\O_{Y,b_i}@>>> L_i,
\end{CD}
$$
where the terms on the right are the fields of fractions. Let $K=k(X)$ and $L=k(Y)$. Then the induced morphism
$$
K_1\otimes_K\ldots\otimes_KK_n\subset L_1\otimes_L\ldots\otimes_LL_n 
$$
is faithfully flat and integral. It follows that each generic point $\eta\in Y_b$ is in the image of $X_a\ra Y_b$.
The latter morphism is integral by Proposition \ref{induced integral}, whence also surjective.
In particular, there is an irreducible closed subscheme $B\subset X_a$ surjecting onto $C\subset Y_b$.
By assumption, $B$ is local henselian, and the morphism $B\ra C$ is surjective and integral.
Whence $C$ is local henselian, according to Proposition \ref{permanence local henselian}.
\qed

\medskip
In particular, an integral normal scheme $X_0$ has the weak Cartan--Artin property if and only
this holds for the total separable closure $X=\TSC(X_0)$. The latter is a scheme that
is everywhere strictly local. For such schemes, the Cartan-- Artin properties
reduce to a striking geometric property with respect to the Zariski topology:

\begin{theorem}
\mylabel{characterization CA}
Let $X$ be an irreducible quasiseparated scheme  that is   everywhere strictly local. 
Then the following are equivalent:
\begin{enumerate}
\item The scheme $X$ has the Cartan--Artin property.
\item The scheme $X$ has the weak Cartan--Artin property.
\item For every pair of irreducible closed subset $A,B\subset X$ 
the intersection $A\cap B$ is irreducible.
\end{enumerate}
\end{theorem}

\proof
Given a sequence $x=(x_1,\ldots,x_n)$  of geometric points, which we may regard as ordinary points $x_i\in X$, we have
  $\O_{X,x_i}=\O_{X,x_i}^s$, whence
$$
X_x=\Spec(\O_{X,x_1})\cap\ldots\cap\Spec(\O_{X,x_n}),
$$
and this is the set of points  in $X$ containing each $x_i$ in their closure.
Since $X$ is quasiseparated, the scheme $X_x$ is quasicompact, so every point
specializes to a closed point. It follows that
$$
X_x=\bigcup_z\Spec(\O_{X,z}),
$$
where the union runs over all closed points $z\in X_x$.

The implication  (i)$\Rightarrow$(ii) is trivial.

To see (iii)$\Rightarrow$(i), it suffices, in light of Theorem \ref{characterization acyclic etale},
to  check that each $X_x$ contains at most one closed point,
that is,  $X_x$ is either empty or local. 
By induction on $n\geq 1$, it is enough to treat the case $n=2$, with $x_1\neq x_2$.
Now suppose there are two closed points $a\neq b$ in $X_x=\Spec(\O_{X,x_1})\cap \Spec(\O_{X,x_2})$, 
and let $A,B\subset X$ be their closure.
Then $x_1,x_2\in A\cap B$. The intersection is irreducible by assumption.
Whence the generic point $\eta\in A\cap B$ lies in $X_x$. Since it is in the closure of  
  $a,b\in X_x$, and these points are closed in $X_x$,  we must have $a=\eta=b$, contradiction.

It remains to verify (ii)$\Rightarrow$(iii). Seeking a contradiction, we suppose that there
are irreducible closed subsets $A,B\subset X$ with $A\cap B$ reducible.
Choose generic points $x_1\neq x_2$ in this intersection, and consider the pair $x=(x_1,x_2)$.
The resulting scheme $X_x=\Spec(\O_{X,x_1})\cap \Spec(\O_{X,x_2})$ is   irreducible, 
because this holds for $X$.
In light of Theorem \ref{characterization acyclic nisnevich}, $X_x$ must be local.
Let $z\in X_x$ be the closed point.
Write $\eta_A,\eta_B$ for the generic points of $A$ and $B$,
respectively. By construction, $\eta_A,\eta_B\in X_x$, whence $z\in A\cap B$. Using that
$x_1,x_2\in A\cap B$ are generic and contained in the closure of $\left\{z\right\}\subset X$, we infer $x_1=z=x_2$, contradiction.
\qed

%===========================================================
\section{Algebraic surfaces}
\mylabel{ALGEBRAIC SURFACES}

Let $k$ be a ground field, and $S_0$ a separated scheme of finite type, assumed to be integral
and 2-dimensional. Note that we do not suppose that $S_0$ is quasiprojective.
Examples of normal surfaces that are proper but not projective appear in \cite{Schroeer 1999}.
In this section we investigate geometric properties of the
absolute integral closure $S=\TSC(S_0)$. Our   result is:

\begin{theorem}
\mylabel{one point}
Assumptions as above. Let $A,B\subset S$ be two integral, 1-dimensional closed subschemes
with $A\neq B$. Then $A\cap B$ contains at most one point.
\end{theorem}

\proof
Clearly, we may assume that $S_0$ is proper and normal, and that the ground field $k$ is algebraically closed.
Employing the notation from   Section \ref{TOTAL SEPARABLE CLOSURE}, we write 
$$
S=\invlim S_\lambda \quadand A=\invlim A_\lambda \quadand B=\invlim B_\lambda, \quad  \lambda\in L.
$$
Seeking a contradiction, we assume that there are two closed points $x\neq y$ in $A\cap B$.
Let $x_\lambda,y_\lambda\in A_\lambda\cap B_\lambda$ be their images on $S_\lambda$.
Enlarging  $\lambda=0$, we may assume that $A_0\neq B_0$ and $x_0\neq y_0$.

The integral Weil divisor $B_0\subset S_0$ has a rational selfintersection number $(B_0)^2\in\QQ$, defined
in the sense of Mumford (\cite{Mumford 1961}, Section II (b)).
Now choose finitely many closed points $z_1,\ldots,z_n\in B_0$ neither contained in $A_0$ nor the singular locus $\Sing(S_0)$,
and let $S'_0\ra S_0$ be the blowing-up with reduced center given by these points.
If we choose $n\gg 0$ large enough, the selfintersection number of the strict transform $B'_0\subset S'_0$ of $B_0$ becomes strictly
negative. According to \cite{Artin 1970}, Corollary 6.12 one may contract it to a point: There exist a proper birational
morphism $g:S_0'\ra Y_0$ with $\O_{Y_0}=g_*(\O_{S'_0})$ sending the integral curve $B'_0$ to a closed point,
and being an open embedding on the complement. Note that $Y_0$ is a proper 2-dimensional algebraic space over $k$, which
is usually not a scheme.
We also write $A_0\subset S'_0$ for the curve $A_0$, considered as a closed subscheme of $S'_0$.
The projection $A_0\ra Y_0$ is a finite   morphism, which is not injective because
the two points $x_0,y_0\in A_0\cap B'_0$ are mapped to a common image. 

We now make an analogous construction in the   inverse system.
Consider the normalization  $S'_\lambda\ra S_\lambda\times_{S_0} S'_0$   of the integral component
dominating $S'_0$. The universal properties of normalization and
fiber product  ensure that the $S'_\lambda$, $\lambda\in L$ form an inverse system, and the transition maps
are obviously affine. In turn, we get an   integral scheme
$$
S'=\invlim S'_\lambda,
$$
coming with a birational morphism $S'\ra S$. This scheme $S'$ is totally separably closed.
Clearly, the projection $S'_\lambda\ra S_\lambda$ are isomorphisms over a neighborhood of the preimage of $A_0\subset S_0$,
which contains $A_\lambda\subset S'$. Therefore,  we may regard
$A_\lambda$ and the points $x_\lambda,y_\lambda$ also as a closed subscheme of $S'_\lambda$.

Let $S'_\lambda\ra Y_\lambda$ be the Stein factorization
of the composition $S'_\lambda\ra S'_0\ra Y_0$. This map contracts the connected components
of the preimage of $B'_0\subset S'_0$ to closed points, and is an open embedding on the complement. 
Using that this preimage contains
the strict transform $B'_\lambda\subset S'_\lambda$
of $B_\lambda\subset S_\lambda$, which is irreducible and thus connected, we infer that
$x_\lambda,y_\lambda\in S'_\lambda$ map to a common image in $Y_\lambda$.
By the universal property of Stein factorizations, the $Y_\lambda$, $\lambda\in L$ form a filtered inverse system.
The transition morphisms $Y_\lambda\ra Y_j$ are finite, because these schemes are finite over $Y_0$. 
In turn, we get another totally separably closed scheme
$Y=\invlim{Y_\lambda}$,
endowed with a birational morphism $S'\ra Y$. 
By construction, the finite morphism $A_\lambda\ra Y_\lambda$ map $x_\lambda,y_\lambda$ to a common image.
Passing to the inverse limit, we get an integral morphism
$$
A=\invlim{A_\lambda}\lra \invlim Y_\lambda=Y,
$$
which maps $x\neq y$ to a common image.
By assumptions,  $A $ is an integral scheme.
According to Proposition \ref{integral => radical}, the map
$A\ra Y$ must be injective, contradiction.
\qed

%===========================================================
\section{Two facts on noetherian schemes}
\mylabel{FACTS ON SCHEMES}

Before we come to higher-dimensional generalizations of Theorem \ref{one point}, we have
to establish two properties pertaining to  quasiprojectivity and connectedness that might
be of independent interest.

\begin{proposition}
\mylabel{chow refined}
Let $R$ be a noetherian ground ring, $X$ a separated  scheme of finite type, and $U\subset X$ be
a quasiprojective open subset. Then there
is a closed subscheme $Z\subset X$ disjoint from $U$ and containing no point $x\in X$ of codimension $\dim(\O_{X,x})=1$
such that the blowing up $f:X'\ra X$ with center $Z$ yields a quasiprojective scheme  $X'$.
\end{proposition}

\proof
By a result of Gross \cite{Gross 2012}, proof for Theorem 1.5, there exists an open subscheme $V\subset X$
containing $U$ and all points of codimension $\leq 1$ that is quasiprojective over $R$.
(Note that his overall assumption that the ground ring $R$ is excellent does not enter in this result.)
By Chow's Lemma (in the refined form of \cite{Deligne 2010}, Corollary 1.4),
there is a blowing-up $X'\ra X$ with center $Z\subset X$ disjoint from $V$ so that
$X'$  becomes quasiprojective over $R$.
\qed

\begin{proposition}
\mylabel{connectedness}
Let $X$ be a   noetherian scheme  satisfying Serre's Condition $(S_2)$, and $D\subset X$ a connected Cartier divisor.
Suppose that the local rings $\O_{X,x}$, $x\in X$ are catenary and homomorphic images of Gorenstein local rings.
Then there is a sequence 
of irreducible components $D_0,\ldots, D_r\subset D$, possibly with repetitions
and covering $D$ such that the successive intersections  $D_{i-1}\cap D_{i}$, $1\leq i\leq r$  have codimension $\leq 2$ in $X$.
\end{proposition}

\proof
Since $D$ is connected and noetherian, there is a a sequence of irreducible components $D_0,\ldots,D_s\subset D$,
possibly with repetitions and covering $D$ so that the successive intersections are nonempty
(\cite{EGA I},  Corollary 2.1.10).  
Fix an index $1\leq i\leq r$, choose closed points $x\in D_{i-1}\cap D_i$, let $R=\O_{X,x}$ be the corresponding local ring on $X$,
and $f\in R$ the regular element defining the Cartier divisor $D$ at $x$.

To finish the proof, we insert between $D_{i-1},D_i$  a sequence of further irreducible components satisfying 
the codimension condition of the assertion on the local ring $R$.
This  can be achieved with Grothendieck's Connectedness Theorem
(\cite{SGA 2}, Expos\'e XIII, Theorem 2.1)  
given in the form of Flenner, O'Carroll and Vogel (\cite{Flenner; O'Carroll; Vogel 1999}, Section 3.1) as follows:

Recall that a noetherian scheme $S$ is called \emph{connected in dimension $d$} if we have $\dim(S)>d$, and
the complement $S\smallsetminus T$
is connected for all closed subsets of dimension $\dim(T)<d$. Set $n=\dim(R)$. Since the local ring $R$ is catenary
and satisfies Serre's condition $(S_2)$,
it must be connected in dimension $n-1$ (loc.\ cit., Corollary 3.1.13).
It follows that $R/fR$ is connected in dimension $n-2$ (loc.\ cit., Lemma 3.1.11, which is essentially 
Grothendieck's Connectedness Theorem).
Consequently, there is a sequence $H_0,\ldots,H_q$ of irreducible components of $\Spec(R/fR)$
so that $H_0,H_q$ are  the local schemes attached to $x \in D_{i-1}, D_i$, and 
whose  successive intersections have dimension $\geq n-2$ in $\Spec(R)$, whence codimension $\leq 2$.
The closures of the $H_0,\ldots,H_q$ inside $D$ yield the desired sequence of further irreducible components.
\qed

\medskip
Let me state this result in a simplified form suitable for most applications:

\begin{corollary}
\mylabel{suitable connectedness}
Let $X$ be a normal irreducible scheme that is of finite type over a regular
noetherian ring $R$, and $D\subset X$ a connected Cartier divisor.
Then there is a sequence 
of irreducible components $D_0,\ldots, D_r\subset D$, possibly with repetitions
and covering $D$ such that the successive intersections  $D_{i-1}\cap D_{i}$, $1\leq i\leq r$  have codimension $\leq 2$ in $X$.
\end{corollary}

\proof
For every affine open subset $U\subset X$, the ring $A=\Gamma(U,\O_X)$ is 
the homomorphic image of some polynomial ring $R[T_1,\ldots,T_n]$. Since the latter
is regular, in particular Gorenstein and catenary, it follows that all local rings
$\O_{X,x}$ are catenary, and  homomorphic images of   Gorenstein local rings.
\qed

%===========================================================
\section{Contractions to points}
\mylabel{CONTRACTIONS}

Let $k$ be a ground field, and $X$ a proper scheme, assumed to be normal and irreducible.
Set $n=\dim(X)$. We say that a connected  closed subscheme $D\subset X$ is \emph{contractible to a point}
if there is a proper algebraic space $Y$ and a  morphism $g:X\ra Y$ with $\O_Y=f_*(\O_X)$ 
sending $D$ to a closed point and being an open embedding on the complement.
Note that such a morphism is  automatically proper, because $X$ is proper and $Y$ is separable.
Moreover, it is unique up to unique isomorphism.

Let $C\subsetneqq X$ be a connected closed subscheme.
The goal of this section is to  construct a projective birational morphism $f:X''\ra X$ so that some
connected Cartier divisor $D''\subset X''$ with $f(D'')=C$ becomes contractible 
to a point, generalizing
the procedure in dimension $n=2$ used in the proof for Proposition \ref{one point}.
The desired morphism $f$ will be a composition  $f=g\circ h$ of two birational morphisms.

The first morphism $g:X'\ra X$ will replace the connected closed subscheme $C\subset X$ by some connected
Cartier divisor $D'\subset X'$ that is a \emph{projective} scheme.
It will be itself a composition of three projective birational morphisms
$$
X \stackrel{g_1}{\longleftarrow} X_1\stackrel{g_2}{\longleftarrow} X_2\stackrel{g_3}{\longleftarrow} X_3=X'.
$$
Here $g_1:X_1\ra X$ is the normalized blowing-up with center $C\subset Y$. Such a map has connected fibers 
by Zariski's Main Theorem. Let $D_1=g_1^{-1}(C)$.
% This is a connected Cartier divisor on $X_1$ by Proposition \ref{}.
According to Proposition \ref{chow refined}, there is a closed subscheme $Z\subset D_1$
containing no point $d\in D_1$ of codimension $\dim(\O_{D_1,d})\leq 1$ so that
the blowing-up $\tilde{D}_1\ra D_1$ becomes projective.  
Let $g_2:X_2\ra X_1$ be the normalized blowing-up with the same center regarded as a closed subscheme $Z\subset X$.
% Note that it has connected fibers, again because $X_1$ is normal.
Viewing $\tilde{D}_1$ as a closed subscheme on the blowing-up $\tilde{X}_1\ra X_1$
with center $Z$, we denote by $D_2\subset X_2$ its preimage on the normalized blowing-up,
which is a   Weil divisor. Finally, let $g_3:X=X_3\ra X_2$ be the normalized blowing-up
with center $D_2\subset X_2$, and $D'=D_3$ be the preimage of $D_2$.
Then $D'\subset X'$ is a Cartier divisor by the universal property of blowing-ups.
The morphism $X'\ra X$ has connected fibers, by Zariski's Main Theorem, and induces
a surjection $D'\ra C$. Since $C$ is connected, it follows that $D'$ is connected.
The proper scheme $D'$ is projective,
because $\tilde{D}$ is projective, and finite maps and blowing-ups are projective morphisms.
Setting $g=g_1\circ g_2\circ g_3:X'\ra X$, we record:

\begin{proposition}
\mylabel{first blowing-up}
The scheme $X'$ is proper and normal, the morphism $g:X'\ra X$ is projective and birational,
the closed subscheme $D'\subset X'$ is a connected Cartier divisor that is a projective 
scheme, and $g(D')=C$.
\end{proposition}

In a second step, we now construct the proper birational morphism $h:X''\ra X'$, so that $X''$
contains the desired contractible effective Cartier divisor.
Choose an ample Cartier divisor $A\subset D'$ containing no generic point from
the intersection of two   irreducible component of $D'$, and so that the invertible
sheaf $\O_{D'}(A-D')$ is ample.
Let $\varphi:\tilde{X}'\ra X'$ be the blowing-up with center $A\subset X'$.
By the universal property of blowing-ups, there is a partial section $\sigma:D'\ra \tilde{X}'$.
Let $\nu:X''\ra\tilde{X}'$ be the normalization map, such that the composition
$$
h:X''\stackrel{\nu}{\lra} \tilde{X}' \stackrel{\varphi}{\lra} X'
$$
is the normalized blowing-up with center $A$. Set $D''=\nu^{-1}(\sigma(D'))\subset X''$,
and let $f=g\circ h:X''\ra X$ be the composite morphism.

\begin{theorem}
\mylabel{second blowing-up}
The scheme $X''$ is proper and normal, the morphism $f:X''\ra X$ is projective and birational,
and $D''\subset X''$ is a connected Cartier divisor that is a projective scheme with $f(D'')=C$.
Moreover, the closed subscheme $D''\subset X''$ is contractible to a point.
\end{theorem}

\proof
By construction, $X''$ is normal and proper, $f$ is projective and birational, and $f(D'')=C$.
According to \cite{Perling; Schroeer 2014}, Lemma 4.4, the subschemes 
$$
\sigma(D'),\varphi^{-1}(A),\varphi^{-1}(D')\subset\tilde{X}'
$$
are Cartier, with $\varphi^{-1}(D)=\sigma(D')+\varphi^{-1}(A)$ in the group of Cartier divisors on
$\tilde{X}'$. Pulling back along the finite surjective morphism between integral schemes $\nu:X''\ra\tilde{X}'$,
we get Cartier divisors $D''=\nu^{-1}(\sigma(D')), E=h^{-1}(A), \hat{D}=h^{-1}(D')$ on $X''$
satisfying $\hat{D}=D''+E$. 
In turn, $\O_{D''}(-D'')\simeq \O_{D''}(E-\hat{D})$. The latter is isomorphic
to the preimage with respect to  the normalization map $\nu:X''\ra\tilde{X}'$ of 
$$
\O_{\sigma(D')}(-1)\otimes \varphi^*\O_{D'}(-D').
$$
The first tensor factor becomes under the canonical identification $\sigma(D')\ra D'$
the invertible sheaf $\O_{D}(A)$. Summing up, $\O_{D''}(-D'')$ is isomorphic to the pull back
of the ample invertible sheaf $\O_{D'}(A-D')$ under a finite morphism, whence is ample.
By Artin's Contractibility Criterion \cite{Artin 1970}, Corollary 6.12, the connected components of $D''$ are
contractible to points.

It remains to verify that $D''$ is connected. For this it suffices to check that the dense
open subscheme  $D''\smallsetminus E$ is connected, which is isomorphic to $D'\smallsetminus A$.
By Proposition \ref{connectedness}, there is a sequence of irreducible components $D'_1,\ldots, D'_r\subset D'$,
possibly with repetitions and covering $D''$ so that the successive intersections $D'_{i-1}\cap D'_i$
have codimension $\leq 2$ in $X'$, whence codimension $\leq 1$ in $D'_{i-1}$ and $D'_i$. 
By the choice of $A\subset D'$, the complements
$(D'_{i-1}\cap D'_i)\smallsetminus A$ are nonempty, and it follows that $D'\smallsetminus A$ is connected.
\qed

%===========================================================
\section{Cyclic systems}
\mylabel{CYCLIC SYSTEMS}

Let $X$ be a scheme. Let us call a pair of closed subschemes
$A,B\subset X$ a \emph{cyclic system} if the following holds:
\begin{enumerate}
\item
The space $A$ is irreducible.
\item
The space $B$ is connected.
\item
The intersection $A\cap B$ is disconnected.
\end{enumerate}

This ad hoc definition will be   useful when dealing with totally separably closed schemes, 
and is somewhat more
flexible than the ``polygons'' used by \cite{Artin 1971}, compare also \cite{Sergio 1985}.
Note that this notion is entirely topological in nature, and
that the conditions ensures that $A,B$ are nonempty and $A\not\subset B$.
One should bear in mind a picture like the following: 
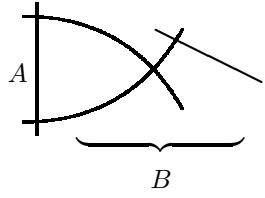
\begin{center}
\begin{picture}(100,80)
\thicklines
\put(25,20){\line(0,1){50}}
\put(14,40){$A$}
\qbezier(20,65)(60,65)(80,30)
\qbezier(20,25)(60,25)(80,60)
\put(70,60){\line(2,-1){40}}
\put(40,20){$\underbrace{\phantom{aaaaaaaaaaaa}}_{}$}
\put(68,0){$B$}
\end{picture}
\end{center}
In this section, we establish some elementary but technical permanence properties for cyclic systems, which
will be used later.

\begin{proposition}
\mylabel{preimage systems}
Let $A,B\subset X$ be a cyclic system, and $f:X'\ra X$ closed surjection with
connected fibers.
Let $A'\subset X'$ be a closed irreducible subset with $f(A')=A$  and $B'=f^{-1}(B)$. Then the pair $A',B'\subset X'$ is 
a cyclic system.
\end{proposition}

\proof
By assumption, $A'$ is irreducible.
Since $f$ is a closed, surjective and continuous, the space $B$ carries the quotient topology
with respect to  $f:B'\ra B$.
Thus $B'$ is connected, because $B$ is connected and $f$ has connected fibers
(\cite{EGA I}, Proposition 2.1.14).
The set-theoretic ``projection formula'' gives
$$
f(A'\cap B')=f(A'\cap f^{-1}(B))=f(A')\cap B=A\cap B,
$$
whence $A'\cap B'$ is disconnected. 
\qed

\begin{proposition}
\mylabel{transform systems}
Let $A,B\subset X$ be a cyclic system, and $f:X'\ra X$ be a     morphism.
Let $U\subset X$ be an open subset over which $f$ becomes an isomorphism, and
let $A',B'\subset X'$ be the   closures of $f^{-1}(A\cap U),f^{-1}(B\cap U)$, respectively.
Suppose the following:
\begin{enumerate}
\item 
$B\cap U$ is connected.
\item
$U$ intersects at least two connected components of $A\cap B$.
\end{enumerate}
Then the pair $A',B'\subset X'$ is  a cyclic system.
\end{proposition}

\proof
Since $A$ is irreducible, so is the nonempty open
subset $A\cap U$, and thus the closure $A'$. Since $B\cap U$ is connected by
Condition (i), the same holds for its closure $B'$.
It remains to check that $A'\cap B'$ is disconnected. 
Since $f$ is continuous and $A\subset X$ is closed, we have 
$$
f(A')=f(\overline{f^{-1}(A\cap U)})\subset\overline{A\cap U}\subset A,
$$
and similarly $f(B')\subset B$. Hence $f(A'\cap B')\subset A\cap B$.
Seeking a contradiction, we assume that $A'\cap B'$ is connected.
Choose points $u,v\in A\cap B$ from two different connected components with $u,v\in U$.
Regarding $u,v$ as elements from $X'$ via $U=f^{-1}(U)$, we see that they are contained
in   the connected subset $A'\cap B'\subset X'$, whence the points $u,v\in A\cap B$ are contained in a connected subset
of $f(A'\cap B')$, contradiction.
\qed

\begin{proposition}
\mylabel{cofinal systems}
Let $X=\invlim X_\lambda$ be a filtered inverse system of quasicompact quasiseparated schemes with integral transition maps,
and  $A,B\subset X$ be a cyclic system.
Then their images $A_\lambda,B_\lambda\subset X_\lambda$ form cyclic systems for
a cofinal subset of indices $\lambda$.
\end{proposition}

\proof
The subsets $A_\lambda,B_\lambda\subset X_\lambda$ are irreducible respective connected, 
because they are images of such spaces.
They are closed, because the projections $X\ra X_\lambda$ are closed maps.
It remains to check that $A_\lambda\cap B_\lambda$ are disconnected for a cofinal subset
of indices. Since inverse limits commute with fiber products, we have
$$
\invlim(A_\lambda)\cap\invlim(B_\lambda)=\invlim(A_\lambda\cap B_\lambda).
$$
The canonical inclusion $A\subset\invlim(A_\lambda)$ of closed subsets
in $X$ is bijective, since the underlying set of the letter is
the inverse limit of the underlying sets of the $X_\lambda$. Similarly we have
$B=\invlim(B_\lambda)$ as closed subsets.
Thus $A\cap B=\invlim(A_\lambda\cap B_\lambda)$.
By assumption, $A\cap B$ is disconnected and the $A_\lambda\cap B_\lambda$ are quasicompact.
It follows from Proposition \ref{decomposition limit} that $A_\lambda\cap B_\lambda$ must
be disconnected for a cofinal subset of indices $\lambda$.
\qed

%===========================================================
\section{Algebraic schemes}
\mylabel{ALGEBRAIC SCHEMES}

We now come to the main result of this paper,
which  generalizes Theorem \ref{one point} from surfaces to arbitrary dimensions. It   takes
the following form:

\begin{theorem}
\mylabel{no cyclic system}
Let $k$ be a ground field and  $X_0$ a integral   scheme that is separated, connected and of finite type. 
Then its total separable closure $X=\TSC(X_0)$ contains no cyclic system.
In particular, for any pair $A,B\subset X$ of irreducible closed subsets,
the intersection $A\cap B$ remains irreducible.
\end{theorem}

\proof
Seeking a contradiction, we assume that there is a cyclic system $A,B\subset X$.
Passing to normalization and compactification, we easily reduce to the situation that $X_0$ is normal and proper.
As usual, write the total separable closure  as an inverse limit of finite $X_0$-schemes $X_\lambda$, $\lambda\in L$
that are connected normal, 
and let $A_\lambda,B_\lambda \subset X_\lambda$  be the
images of $A,B\subset X$, respectively. Then we have
$$
X=\invlim X_\lambda,\quad A=\invlim A_\lambda\quadand B=\invlim B_\lambda,
$$
and all transition morphisms 
are finite and surjective.  
Clearly, the $A_\lambda$ are irreducible and the $B_\lambda$ are connected. In light of Proposition \ref{cofinal systems},
we may replace $L$ by a cofinal subset and assume that the $A_\lambda,B_\lambda\subset X_\lambda$ form cyclic systems
for all $\lambda\in L$, including $\lambda=0$.

To start with, we deal with a rather \emph{special case}, namely that the connected subset 
$B_0\subset X_0$ is \emph{contractible to a point}. Let $X_0\ra Y_0$ be the contraction.
Then the composite map $A_0\ra Y_0$ has a disconnected fiber, because the connected
components of the disconnected subset $A_0\cap B_0$ are mapped to a common image point.
Next, consider the induced maps $X_\lambda\ra Y_\lambda$ contracting the connected components
of the preimage $X_\lambda\times_{X_0}B_0\subset X_\lambda$ to points. Again the composite map $A_\lambda\ra Y_\lambda$ has a disconnected fiber. 
This is because the irreducible subset $A_\lambda$ is not
contained in the preimage $X_\lambda\times_{X_0}B_0$,   the connected subset $B_\lambda$ is contained in $X_\lambda\times_{X_0}B_0$, 
thus contained in a connected component of the latter, whence maps to a point in the scheme $Y_\lambda$.

To proceed, consider the Stein factorization $A_\lambda\ra A_\lambda'\ra Y_\lambda$. Then $A_\lambda'\ra Y_\lambda$
is finite, but not injective.
Passing to the inverse limits $A'=\invlim A_\lambda'$ and $Y=\invlim Y_\lambda$,
we obtain a totally separably closed scheme $Y$,  an integral scheme $A'$,
and an integral morphism $A'\ra Y$.
The latter must be radical, by Proposition \ref{integral => radical}.
In turn, the maps $A_\lambda'\ra Y_\lambda$ are radical for sufficiently large indices $\lambda\in L$,
according to \cite{EGA IVc}, Theorem 8.10.5. In particular, these maps are injective,  contradiction.

We now come to the \emph{general case}.  
We shall proceed in six steps, some of which are repeated. In each step,  a given proper normal connected scheme $X_0$
whose absolute separable closure $X=\TSC(X_0)$ contains a cyclic system $A,B\subset X$
will be replaced by another such scheme with slightly modified properties, until we finally 
reach the special case discussed in the preceding paragraphs. As we saw, the latter is impossible.
To start with,   choose points $u,v\in A\cap B$ coming from different connected components.
Let $u_\lambda,v_\lambda\in A_\lambda\cap B_\lambda$   be their images, respectively.

\medskip
{\bf Step 1:}
\emph{We first reduce to the case that   $u_0,v_0\in A_0\cap B_0$ lie in different connected components.}
Suppose that for all indices $\lambda\in L$, the points $u_\lambda,v_\lambda\in A_\lambda\cap B_\lambda$ are contained
in the same connected component $C_\lambda\subset A_\lambda\cap B_\lambda$.
Then the $C_\lambda$, $\lambda\in L$ form an inverse system with finite transition maps,
and their inverse limit $C=\invlim C_\lambda$ is connected 
(\cite{EGA IVc}, Proposition 8.4.1, compare also Appendix A).
Using $u,v\in C\subset A\cap B$, we reach a contradiction.
Hence there exist an index $\lambda\in L$ such that $u_\lambda,v_\lambda\in A_\lambda\cap B_\lambda$ stem from different
connected components. 
Replacing $X_0$ by such an $X_\lambda$ we thus may assume that  $u_0,v_0\in A_0\cap B_0$ lie
in different connected components.

\medskip
{\bf Step 2:} 
\emph{ Next we reduce to the case that $u_0,v_0\in A_0\cap B_0$ are not contained in a common
irreducible component of $B_0$.}
Let $f:X'_0\ra X_0$ be the normalized blowing-up with center $A_0\cap B_0\subset X_0$, such that
the strict transforms of $A_0,B_0$ on $X'_0$ become disjoint. 

Let $X'_\lambda$ be the normalization
of $X_\lambda\times_{X_0}X'_0$ inside the function field $k(X_\lambda)=k(X'_\lambda)$.
Then the $X'_\lambda$, $\lambda\in L$ form an inverse system of proper normal connected schemes
whose transition morphisms are finite and dominant. Let $X'$ be their inverse limit.
By construction, we have projective birational morphisms $X'_\lambda\ra X_\lambda$ and
a resulting   birational morphism $X'\ra X$, which is a proper surjective map  with connected fibers,
according to \cite{EGA IVc}, Proposition 8.10.5 (xii), (vi) and (vii). 
Note that all these $X$-morphisms become isomorphisms when base-changed to $X_0\smallsetminus (A_0\cap B_0)$.

Let $A' \subset X'$ be the strict transform of $A\subset X$, and $B'\subset X'$
the reduced scheme-theoretic preimage of $B\subset X$.
These form a cyclic system on $X'$, by Proposition \ref{preimage systems}.
Consider their images $A'_\lambda,B'_\lambda\subset X'_\lambda$, which are closed
subsets. The $A'_\lambda$ are irreducible and the $B'_\lambda$ are connected.
Clearly, the $A'_\lambda$ can be viewed as the strict transforms of $A_\lambda\subset X_\lambda$.
Moreover, the inclusion $B'_\lambda\subset X'_\lambda\times_{X_\lambda}B_\lambda$ is an equality of sets. This is because
the morphisms  $B\ra B_\lambda$  and $X'\ra X$ are surjective,
by \cite{EGA IVc}, Theorem 8.10.5 (vi).
Again with Proposition \ref{preimage systems}, we infer that all $A_\lambda,B_\lambda\subset X_\lambda$ form
a cyclic system.

Let $\varphi:X'\ra X$ be the canonical morphism. The set-theoretical Projection Formula
gives $\varphi(A'\cap B')=A\cap B$. Choose points $u',v'\in A'\cap B'$ mapping
to $u$ and $v$, respectively, and consider their images $u'_0,v'_0\in A'_0\cap B'_0$.
The latter are not contained in the same connected component, because they map to $u_0,v_0\in A_0\cap B_0$.
Moreover, we claim that the $u'_0,v'_0\in A'_0\cap B'_0$ are not contained in a common irreducible component of $B'_0$:
Suppose they would lie in a common irreducible component $W_0'\subset B'_0$.
Since the strict transforms of $A_0,B_0$
on $X'_0$ are disjoint, the image $f(W'_0)\subset X_0$ under the blowing-up $f:X'_0\ra X_0$ must lie in the center
$A_0\cap B_0$. Then $u_0,v_0\in f(W_0)\subset A_0\cap B_0$ are contained in some connected subset, contradiction.

Summing up, after replacing $X_0$ by $X'_0$, we may additionally assume that  $u_0,v_0\in B_0$ 
are not contained in a common irreducible component of $B_0$.

\medskip\noindent
{\bf Step 3:}
\emph{The next goal is to turn  the closed subscheme  $B_0\subset X_0$ into a Cartier divisor.}
Let $X'_0\ra X_0$ be the normalized blowing-up with center $B_0\subset X_0$, such that
its preimage on $X'_0$ becomes a Cartier divisor. 
Define $X'_\lambda, X',A',B',A'_\lambda,B'_\lambda$ as in Step 2.
Then $A',B'\subset X'$ and the $A'_\lambda,B'_\lambda\subset X'_\lambda$ form cyclic systems.

Again let $\varphi:X'\ra X$ be the canonical morphism. Then
$\varphi(A'\cap B')=A\cap B$ by the projection formula. Let $u',v'\in A'\cap B'$
be in the preimage of the points $u,v\in A\cap B$.
Since their images $u_0,v_0\in A_0\cap B_0$ do not lie in a common connected component and are not contained
in a common irreducible component of $B$,
the same necessarily holds for their images $u'_0,v'_0\in A'_0\cap B'_0$.

Replacing $X_0$ by $X'_0$, we thus may additionally assume that $B_0\subset X_0$ is the support of
an effective Cartier divisor.

\medskip
{\bf Step 4:} 
\emph{We now reduce to the case that all $B_\lambda\subset X_\lambda$ are supports of connected effective Cartier divisors.}
The preimages $X_\lambda\times_{X_0}{B_0}\subset X_\lambda$ of the Cartier divisor $B_0\subset X_0$ remain Cartier divisors,
because $X_\lambda\ra X$ is a dominant morphism between normal schemes.
Let $C_\lambda\subset X_\lambda\times_{X_0}B_0$ be the connected component containing the connected subscheme 
$B_\lambda\subset X_\lambda\times_{X_0}B_0$.
Clearly, the $C_\lambda\subset X_\lambda$ are connected Cartier divisors, which form an inverse system such
that $C=\invlim C_\lambda$ contains $B=\invlim B_\lambda$. Note that the transition maps $C_j\ra C_\lambda$, $i\leq j$
and the projections $C\ra C_\lambda$ are not necessarily surjective. 
Nevertheless, we conclude with \cite{EGA IVc}, Proposition 8.4.1 that
$C$ is connected. Clearly, the points $u,v\in A\cap C$  
lie in different connected components, and are not contained in a common irreducible component of $C$,
because this hold for their images $u_0,v_0\in A_0\cap B_0$. Thus $A,C\subset X$  form
a cyclic system.

Replacing $B$ by $C$ and $B_\lambda$ by $C_\lambda$, we thus may
additionally assume that all $B_\lambda\subset X_\lambda$ are supports of connected Cartier divisors.
Note that at this stage we have sacrificed our initial assumption that the the  projections $B\ra B_\lambda$ 
and the transition maps $B_j\ra B_\lambda$
are surjective.

\medskip
{\bf Step 5:}
\emph{Here we reduce to the case that the proper scheme $B_0$ is projective.}
Since $u_0,v_0\in B_0$ are not contained in a common irreducible component, they
admit a common affine open neighborhood inside $B_0$: To see this, let
$B_0'\subset B_0$ be the finite union of the irreducible components that do not contain $u_0$,
and let $B_0''\subset B_0$ be the finite union of the irreducible components that do not contain $v_0$.
Then $B_0=B_0'\cup B_0''$, so any affine open neighborhoods of $u_0\in B_0\smallsetminus B_0'$ and $v_0\in B_0\smallsetminus B_0''$
are disjoint, and their union is the desired common affine open neighborhood.

It follows that there is a closed subscheme $Z_0\subset B_0$
containing neither $u_0,v_0$ nor any point $b\in B_0$ of codimension $\dim(\O_{B_0,b})=1$
so that the blowing-up $\tilde{B}_0\ra B_0$ with center $Z_0$
yields a projective scheme $\tilde{B}_0$. This holds by Proposition \ref{chow refined}.
Note that $Z_0\subset B_0$ has codimension $\geq 2$.

Let $X'_0\ra X_0$ be the normalized blowing-up with the same center $Z_0$, regarded as a closed subscheme $Z_0\subset X_0$.
Define $X'_\lambda$ and $X'=\invlim X'_\lambda$ as in Step 2.  
Let $U\subset X$ be the complement of the closed subscheme $X\times_{X_0}Z_0\subset X$.
This open subscheme intersects at least two connected components of $A\cap B$, because $u,v\in U$.
Clearly, the canonical morphism $f:X'\ra X$ becomes an isomorphism over $U$.
Furthermore, the intersection $B\cap U$ is connected. To see this, set $Z_\lambda=X_\lambda\times_{X_0}Z_0$
and $U_\lambda=X_\lambda\smallsetminus Z_\lambda$, such that $B\cap U=\invlim(B_\lambda\cap U_\lambda)$.
In light of \cite{EGA IVc}, Proposition 8.4.1, it suffices to check that $B_\lambda\cap U_\lambda$ is
connected for each $\lambda\in L$. Indeed, by the Going-Down Theorem applied to the finite morphism
$B_\lambda\ra B_0$, all $Z_\lambda\subset B_\lambda$ have
codimension $\geq 2$. The $B_\lambda\subset X_\lambda$ are supports of connected Cartier divisors.
In light of Corollary \ref{suitable connectedness}, 
the set $B_\lambda\cap U_\lambda=B_\lambda\smallsetminus Z_\lambda$ must remain connected.

According to Proposition \ref{transform systems}, the closures $A',B'\subset X'$  of $A\cap U,B\cap U$ form
a cyclic system. Clearly, the image of the composite map $B'\ra X_0$ is $B_0\subset X_0$.
Replacing $X$ by $X'$, we  may   assume that $B_0$  is a projective scheme.
Note that with this step we have sacrificed the property that the $B_\lambda\subset X_\lambda$ are
Cartier divisors.

\medskip
{\bf Step 6:} 
\emph{We now repeat Step  4, and achieve again that $B_\lambda\subset X_\lambda$ are supports of connected Cartier divisors.} 
Moreover, the $B_\lambda$ stay projective, because $B_0$ is projective and the $X_\lambda\ra X_0$ are finite.

\medskip
{\bf Step 7:} 
\emph{In this last step we achieve that   $B_0\subset X_0$ becomes contractible to a point, thus reaching a contradiction.}
Choose an ample Cartier divisor $H_0\subset B_0$ containing no generic point
from the intersection of two irreducible components of $B_0$,  and no generic point from
the intersection $A_0\cap B_0$, and no point from $\left\{u_0,v_0\right\}$. Furthermore, we demand that 
the invertible $\O_{B_0}$-module $\O_{B_0}(H_0-B_0)$ is ample. Let $X'_0\ra X_0$ be the normalized  blowing-up
with center $H_0\subset X_0$, and $A'_0,B'_0\subset X'_0$ the strict transform of $A_0,B_0\subset X_0$, respectively.
As in the proof for Theorem \ref{second blowing-up}, the subscheme $B'_0\subset X'_0$ is connected and contractible to a point.
Define $X'_\lambda$ and $X'=\invlim X'_\lambda$ as in Step 2. 

Let $U\subset X$ be the complement of the closed subscheme $X\times_{X_0}H_0$.
This open subscheme intersects at least two connected components of $A\cap B$, because $u,v\in U\cap V$.
Clearly, the canonical morphism $f:X'\ra X$ becomes an isomorphism over $U$.
Furthermore, the intersection $B\cap U$ is connected. To see this, set $H_\lambda=X_\lambda\times_{X_0}H_0$
and $U_\lambda=X_\lambda\smallsetminus H_\lambda$, such that $B\cap U=\invlim(B_\lambda\cap U_\lambda)$.
In light of \cite{EGA IVc}, Proposition 8.4.1, it suffices to check that $B_\lambda\cap U_\lambda=B_\lambda\smallsetminus H_\lambda$ is
connected for each $\lambda\in L$. 
Since $B_\lambda\subset X_\lambda$ is the support of a connected Cartier divisor, there is a sequence
of irreducible components $C_{1},\ldots,C_{n}$ covering $B_\lambda$ so that
each successive intersection $C_{j-1}\cap C_{j}$ has codimension $\leq 2$ in $X_\lambda$.
Let $f:X_\lambda\ra X_0$ be the canonical   morphisms, and consider the inclusion
$$
f(C_{j-1}\cap C_{j})\subset f(C_{j-1})\cap f(C_j).
$$
Since $f$ is finite and dominant, the left hand side has codimension $\leq 2$ in $X_0$,
and the $f(C_{j-1}), f(C_j)$ are irreducible components from $B_0$.
By the choice of the ample Cartier divisor $H_0\subset X_0$, 
the image  $f(C_{j-1}\cap C_{j})$
is not entirely contained in $H_0$. We conclude that $(C_{j-1}\cap C_{j})\smallsetminus H_\lambda$ is 
nonempty, thus $B_\lambda\smallsetminus H_\lambda$ remains connected.

Let $A',B'\subset X'$ be the closure of $A\cap U,B\cap U$.
The former comprise a cyclic system, in light of Proposition \ref{transform systems}.
Clearly, the image of $B'\ra X_0$ equals $B'_0$. Replacing $X_0$ by $X_0'$,
we thus have achieved the situation that $B_0\subset X_0$ is contractible to a point.
We have seen in the beginning of this proof that such a cyclic system $A,B\subset X$
does not exist.                                                                                         
\qed

%===========================================================
\section{Nisnevich \v{C}ech cohomology}
\mylabel{NISNEVICH CECH}

We now apply our results to prove the following:

\begin{theorem}
\mylabel{cech = sheaf cohomology}
Let $k$ be a ground field, and  $X$ be a quasicompact and separated $k$-scheme.
Then $\cH^p(X_\Nis,\underline{H}^q(F))=0$ for all integers $p\geq 0$, $q\geq 1$  and
all abelian Nisnevich sheaves $F$ on $X$.
In particular, the canonical maps
$$
\cH^p(X_\Nis,F)\lra H^p(X_\Nis,F) 
$$
are bijective for all $p\geq 0$.
\end{theorem}

\proof
Recall that 
$\check{H}^p(X_\Nis, \underline{H}^q(F)) = \dirlim_{U} \check{H}^p(U_\Nis,\underline{H}^q(F))$,
where the direct limit runs over the system of all completely decomposed \'etale surjections $U\ra X$.
For an explanation of the transition maps, we refer to \cite{Godement 1964}, Chapter II, Section 5.7.
Since $X$ is quasicompact and \'etale morphisms are open, one easily sees that
the the subsystem of all \'etale surjections $U\ra X$ with $U$ affine form a cofinal subsystem.

Let $[\alpha]\in \check{H}^p(X_\Nis, \underline{H}^q(F))$ be a \v{C}ech class with $q\geq 1$.
Represent the class by some cocycle $\alpha\in H^q(U^p,F)$, where $U\ra X$ is a completely decomposed
\'etale surjection,
and $U^p$ denotes the $p$-fold selfproduct in the category of $X$-schemes.
This selfproduct   remains quasicompact, because $X$ is quasiseparated.
Since $q\geq 1$, there is a completely decomposed \'etale surjection $W\ra U^p$ with $W$ affine so that $\alpha|W=0$.
It suffices to show that there is a refinement $U'\ra U$ so that 
the structure morphisms $U'^p\ra U^p$ factors over $W\ra U^p$, because then the class of $\alpha$
in the direct limit of the \v{C}ech complex vanishes. From this point on we merely work
with the schemes $W,U,X$  and may forget about
the   sheaf $F$ and the cocycle $\alpha$.

This allows us to reduce to the case of $k$-schemes of finite type:
Write $X=\invlim X_\lambda$ as an inverse limit with affine transition maps of 
$k$-schemes $X_\lambda$ that are of finite type.
Passing to a cofinal subset of indices, we may assume that there is a smallest index $\lambda=0$, so that
$X_0$ is separated (\cite{Thomason; Trobaugh 1990}, Appendix C, Proposition 7).
Since $U$ is quasicompact,  our morphisms $U\ra X$ is of finite presentation,
so we  may assume that $U=U_0\times _{X_0}X_0$ for some $X_0$-scheme $U_0$ of finite presentation.
Similarly, we can assume $W=W_0\times_{U_0^p}U^p$ for some $U_0^p$-scheme $W_0$ of finite presentation.
We may assume that $U_0\ra X_0$ and $W_0\ra U_0^p$ are \'etale (\cite{EGA IVd}, Proposition 17.7.8) and surjective
(\cite{EGA IVc}, Proposition 8.10.5).
According to Lemma \ref{limit completely decomposed} below, we can also impose that these morphisms are completely decomposed.
Summing up, it suffices to produce  a refinement $U_0'\ra U_0$ so that $U_0'^p\ra U_0^p$ factors over $W_0$.
In other words, we may assume that $X$ is of finite type over the ring $\ZZ$.
In particular, $X$  has only finitely many irreducible components,   the normalization
of the corresponding integral schemes are finite over $X$, and $X$ is noetherian.

Next, we make use of the weak Cartan--Artin property.
Indeed, according to Artin's induction argument in  \cite{Artin 1971}, Theorem 4.1,
applied to the Nisnevich topology rather than the \'etale topology,
we see that it suffices to check that $X$ satisfies the 
weak Cartan--Artin property. Let $X_0$ be the normalization of some irreducible component
of $X$, viewed as an integral scheme.
According to Propositions \ref{irreducible Cartan--Artin} and \ref{normalization Cartan--Artin}, it is permissible to replace $X$ by $X_0$.
Changing notation, we write $X=\TSC(X_0)$ for the total separable closure of the integral scheme $X_0$.
In light of Theorem \ref{image weak cartan--artin}, it is enough to verify the weak Cartan--Artin property for $X$.

Here, our problem reduces to a simple geometric statement: According to Theorem \ref{characterization CA},
it suffices to check that for each pair of irreducible closed subsets $A,B\subset X$, the 
intersection $A\cap B$ remains irreducible. 
And indeed, this holds by Theorem \ref{no cyclic system}. Note that only in this very last step,
we have used that $X$ is separated, and  the existence of a ground field $k$.
\qed

\medskip
In the preceding proof, we have used the following facts:

\begin{lemma}
\mylabel{limit completely decomposed}
Suppose that $X=\invlim X_\lambda$ is a filtered inverse system of quasicompact schemes with affine transition maps,
and $U\ra X$ is a completely decomposed \'etale surjection, with  $U$   quasicompact. Then  there is an 
index $\alpha$ and a completely decomposed \'etale surjection $U_\alpha\ra X_\alpha$
with $U=X\times_{X_\alpha}U_\alpha$.
\end{lemma}

\proof
The morphism $U\ra X$ is of finite presentation, because it is \'etale and its domain is
quasicompact. According to \cite{EGA IVc}, Theorem 8.8.2 there is an $X_\alpha$-scheme $U_\alpha$ of finite presentation
with $U=X\times_{X_\alpha}U_\alpha$ for some index $\alpha$.  We may assume that $U_\alpha\ra X_\alpha$ is surjective
(\cite{EGA IVc}, Theorem 8.10.5) and \'etale (\cite{EGA IVd}, Proposition 17.7.8).

For each index $\lambda\geq \alpha$, let $F_\lambda\subset X_\lambda$ be the subset of all points $s\in X_\lambda$ so that
the \'etale surjection $U_\lambda\otimes\kappa(s)\ra\Spec\kappa(s)$ is completely decomposed, that is,
admits a section. According to Lemma \ref{ind-constructible} below, this subset is ind-constructible.
In light of \cite{EGA IVd}, Corollary 8.3.4, it remains to show that
the inclusion $\bigcup_\lambda u_\lambda^{-1}(F_\lambda)\subset X$ is an equality,
where $u_\lambda:X\ra X_\lambda$ denote the projections.

Fix a point $x\in X$, and choose a point $u\in U$ above so that the inclusion $k(x)\subset \kappa(u)$
is an equality. Let $x_\lambda\in X_\lambda$ and $u_\lambda\in U_\lambda$ be the respective image points.
Then $U_x=\invlim (U_\lambda)_{x_\lambda}$, because inverse limits commute
with fiber products. Furthermore, the section of $U_x\ra\Spec\kappa(x)$
corresponding to $u$ comes from a section of $(U_\lambda)_{x_\lambda}\ra\Spec\kappa(x_\lambda)$,
for some index $\lambda$, by \cite{EGA IVc}, Theorem 8.8.2.
\qed

\medskip
Recall that a subset $F\subset X$ is called \emph{ind-constructible} 
if each point $x\in X$ admits an open neighborhood $x\in U$ so that $F\cap U$
is the   union of  constructible subsets of $U$
(compare \cite{EGA I}, Chapter I, Section 7.2 and Chapter 0, Section 2.2).

\begin{lemma}
\mylabel{ind-constructible}
Let $f:U\ra X$ be a     \'etale surjection. Then the subset $F\subset X$
of points $x\in X$ for which $U_x\ra\Spec\kappa(x)$ admits a section is ind-constructible.
\end{lemma}

\proof
The problem is local, so we may assume that $X=\Spec(R)$ is affine.
Fix a point $x\in X$ so that $U_x\ra\Spec\kappa(x)$ admits a section, and let $A\subset X$ be its closure.
The section extends to a section for $U_A\ra A$ over some dense open subset $C\subset A$.
Shrinking $C$ if necessary, we may assume that $C=A\cap\Spec(R_g)$ for some $g\in R$.
Then $C\subset X$ is constructible and contains $x$.
\qed

%===========================================================
\section{Appendix A: Connected components of  schemes}
\mylabel{CONNECTED COMPONENTS}

Let $X$ be a topological space and $a\in X$ a point.
The corresponding \emph{connected component} $C=C_a$  is   the union of
all connected subsets containing $a$. This is then the largest connected subset
containing $a$, and it is closed but not necessarily open.
In contrast, the \emph{quasicomponent} $Q=Q_a$   is defined as the intersection
of all  open-and-closed neighborhoods of $a$, which is also closed but
not necessarily open.
Clearly, we have $C\subset Q$. Equality holds if $X$ has only finitely many quasicomponents,
or if $X$ is compact or locally connected, but in general there is a strict inclusion.

The simplest example with $C\subsetneqq Q$ seems to be the following:
Let $U$ be an infinite discrete space, $Y_i=U\cup\left\{a_i\right\}$ be two copies of
the Alexandroff compactification, and $X=Y_1\cup Y_2$ the space obtained by gluing
along the open subset $U$, which is quasicompact but not Hausdorff. Then all connected components are singletons,
but the discrete space $Q=\left\{a_1,a_2\right\}$ is a quasicomponent.

Assume that $X$ is a locally ringed space. Write $R=\Gamma(X,\O_X)$
for the ring of global sections and $\mathfrak{p}\subset R$ for the prime ideal
of all global sections vanishing at our point $a\in X$.
Obviously, the open-and-closed neighborhoods $U\subset X$  
correspond to the idempotent elements $\epsilon\in R\smallsetminus\mathfrak{p}$,
via $U=X_\epsilon=\left\{x\in X\mid\epsilon(x)=1\right\}$.
Let $L\subset R\smallsetminus\mathfrak{p}$ the multiplicative subsystem of all idempotent
elements. We put an order relation
by declaring $\epsilon\leq\epsilon'$ if $\epsilon'\mid\epsilon$. 
Then $\epsilon\mapsto X_\epsilon$ is a filtered inverse system of subspaces,
and  the quasicomponent is 
$$
Q=\bigcap_{\epsilon\in L} X_\epsilon=\invlim_{\epsilon\in L} X_\epsilon.
$$
Now suppose that $X$ is a scheme. Then each inclusion morphism
$\iota_{\epsilon}:X_\epsilon\ra X$ is affine, hence $\shA_\iota=\iota_{\epsilon *}(\O_{X_\epsilon})$
are quasicoherent, and the same holds for $\shA=\dirlim\shA_\epsilon$. 
In turn, we have
$Q=\invlim_{\epsilon\in L} X_\epsilon = \Spec(\shA)$,
which puts a scheme structure on the  quasicomponent. The quasicoherent ideal sheaf
for the closed subscheme $Q\subset X$ is   $\shI=\bigcup_{\epsilon\in L}(1-\epsilon)\O_X$.
The following observation is due to Ferrand in the affine case, and Lazard in the general case
(see  \cite{Lazard 1967}, Proposition 6.1 and Corollary 8.5): 

\begin{lemma}
\mylabel{component = quasicomponent}
Let $X$ be a quasicompact and quasiseparated scheme and $a\in X$ be a point.
Then we have an equality $C=Q$ between the   connected component
and the quasicomponent containing $a\in X$.
\end{lemma}

If $X=\Spec(R)$ is affine, and $a\in X$ is the  point corresponding to a prime ideal
$\mathfrak{p}\subset R$, then the corresponding $C=Q$ is defined by the 
ideal $\mathfrak{a}=\bigcup Re$, where $e$ runs through the idempotent elements in $\mathfrak{p}$.
Clearly, we have $R/\mathfrak{a}=\dirlim_eA/eA$.
The latter description was already given by Artin in  \cite{Artin 1971}, page 292.

Let $X_0$ be a  scheme, and $X=\invlim_{\lambda\in L} X_\lambda$ be an inverse system of affine $X_0$-schemes.
The argument for the previous Lemma essentially depends on 
\cite{EGA IVc}, Proposition 8.4.1 (ii).
There, however, it is incorrectly claimed that a sum decomposition of $X$ comes 
from a sum decomposition of $X_\lambda$
for sufficiently large $\lambda\in L$, provided that $X_0$ is merely quasicompact.
This only becomes true  only under the additional hypothesis of quasiseparatedness:

\begin{proposition}
\mylabel{decomposition limit}
Let $X_0$ be a quasicompact and quasiseparated scheme, and $X_\lambda$ an filtered inverse system
of affine $X_0$-schemes, and $X=\invlim X_\lambda$. If $X=X'\cup X''$ is a decomposition into
disjoint open subsets, that there is some $\lambda\in L$ and a decomposition $X_\lambda=X_\lambda'\cup X''_\lambda$
into disjoint open subset so that $X',X''\subset X$ are the respective preimages.
\end{proposition}

\proof
This is a special case of \cite{EGA IVc}, Theorem 8.3.11, but I would
like to give an alternative proof relying on absolute noetherian approximation
rather than  ind- and pro-constructible sets:
Choose  a $\ZZ$-scheme of finite type $S$ so that there is an affine morphism $X_0\ra S$
(\cite{Thomason; Trobaugh 1990}, Appendix C, Theorem 9).
Replacing $X_0$ by $S$, we may assume that $X_0$ is noetherian.
Let $\shA_\lambda$ be the quasicoherent $\O_{X_0}$-algebras
corresponding to the $X_\lambda$, such that $X$ is the relative spectrum of
$\shA=\dirlim(\shA_\lambda)$.
Since the topological space of $X_0$ is noetherian,  the canonical map
$$
\dirlim H^0(X_0,\shA_\shL)\lra H^0(X_0,\dirlim\shA_\lambda)
$$
is bijective, see \cite{Hartshorne 1977}, Chapter III, Proposition 2.9.  Thus each idempotent $e\in\Gamma(X_0,\shA)$ 
comes from an element $e_\lambda\in \Gamma(X_0,\shA_\lambda)$. Passing to a larger index, we may assume that $e_\lambda^2=e_\lambda$,
such that our element $e_\lambda$ becomes idempotent.
The statement follows from  the aforementioned 
correspondence between idempotent global sections and open-and-closed subsets.
\qed

\medskip
Here is a counterexamples with $X_0$   not quasiseparated: Let $U$ be an infinite discrete set, and $X=Y_1\cup Y_2$ be the
gluing of two copies $Y=Y_1=Y_2$ of the Alexandroff compactification along the open subset $U\subset Y_i$
mentioned above.
We may regard the Alexandroff compactification 
as a profinite space, by choosing a total order on $U$: Let $L$ be the set of all finite subsets of $U$, ordered by inclusion.
Given $\lambda\in L$, we write $F_\lambda\subset U$ for the corresponding finite subset.
For every $\lambda\leq\mu$, let  $F_\mu\ra F_\lambda$ by the retraction of the canonical inclusion
$F_\lambda\subset F_\mu$ sending the complementary points   to the largest element   $f_\mu\in F_\mu$.
Then one easily sees that $Y=\invlim F_\lambda$, where  the point   at infinity becomes the
tuple of largest elements $a=(f_\lambda)_{\lambda\in L}$.

To proceed, fix  a ground field $k$. Then we may endow the profinite space $Y$ with the structure of an affine $k$-scheme,
by regarding   $F_\lambda$ as   the spectrum of the $k$-algebra $R_\lambda=\Hom_\text{Set}(F_\lambda,k)$.
In turn, we regard the gluing $X=Y_1\cup Y_2$ as a $k$-scheme. This scheme is quasicompact, but
not quasiseparated, because the intersection $Y_1\cap Y_2=U$ is   infinite and discrete.
Consider the   subset $V_\lambda =X\smallsetminus F_\lambda$, which is open-and-closed.
Moreover, the inclusion morphism $V_\lambda\ra X$ is affine, and   $\invlim_{\lambda\in L} V_\lambda=Q=\left\{a_1,a_2\right\}$
is the disconnected quasicomponent. The sum decomposition of the quasicomponent $Q$, however, does not come from a sum decomposition 
of any $V_\lambda$.

Here is a counterexample with  $X_0$ not quasicompact:
Let $X_0$ be the disjoint union of two infinite chains  $C'=\bigcup C'_n$ and $C''=\bigcup C''_n$, $n\in\ZZ$
of copies of the   projective line over the ground field $k$, together with further copies $B_n$ of the projective line
joining the intersection $C'_{n-1}\cap C'_n$ with $C''_{n-1}\cap C''_n$.
Let $L$ be the collection of all finite subset of $\ZZ$, ordered by inclusion.
Consider the inverse system of closed subsets $X_\lambda\subset X$ 
consisting of $C'\cup C''$ and the union $\bigcup_{n\not\in\lambda}B_n$. Then all $X_\lambda$ are
connected, but $X=\invlim(X_\lambda)=\bigcap X_\lambda=C'\amalg C''$ becomes disconnected.

%===========================================================
\section{Appendix B: \v{C}ech   and sheaf cohomology}
\mylabel{CECH AND SHEAF}

Here we briefly recall H.\ Cartan's Criterion for equality for \v{C}ech and sheaf cohomology,
which is   described in Godement's monograph \cite{Godement 1964}, Chapter II, Theorem 5.9.2.
Throughout we  work in the context of topoi and sites. The general reference is \cite{SGA 4a}.
Suppose $\shE$ is a topos.
Choose a site $\shC$  and an equivalence $\operatorname{Sh}(\shC)=\shE$, and make the choice 
so that the category $\shC$ has a terminal object $X\in\shC$ and all fiber products, and so
that the  Grothendieck topology is given in terms of a pretopology of coverings.
Suppose we additionally have a   subcategory $\shB\subset\shC$ such that every object
$U\in \shC$ admits a covering $(U_\alpha\ra U)_\alpha$ with $U_\alpha\in \shB$, and that $\shB$ has all fiber products.
Then the subcategory inherits a pretopology, and the restriction functor of sheaves on $\shC$ to 
sheaves on $\shB$ is an equivalence.

For each abelian sheaf $F$ on  $\shC$  and each object $U\in\shC$,
we write $H^p(U,F)$ for the cohomology groups in the sense of derived functors,
and 
$$
\cH^p(U,F)=\dirlim \check{H}^p(\mathfrak{V},\check{H}^q(F)),
$$ 
for the \v{C}ech cohomology groups, where the direct limit runs over
all coverings   $\mathfrak{V}= (U_\alpha\ra U)_\alpha$, and the transition maps are defined 
as in \cite{Godement 1964}, Chapter II, Section 5.7. 
Note that me may restrict to those coverings with all $U_\alpha\in \shB$.

Let $F$ be an abelian sheaf on the site $\shC$.
We denote by $\underline{H}^q(F)$ the presheaf $U\mapsto H^q(U,F)$.
The functors $F\mapsto \underline{H}^p(F)$, $p\geq 0$ form a   $\delta$-functor
from the category of abelian sheaves to the category of abelian presheaves.
This $\delta$-functor is universal, because it vanishes on injective objects for $p\geq 1$.

The  inclusion  functor $F\mapsto\underline{H}^0(F)$
from the category of abelian sheaves to the category of abelian presheaves is 
 right adjoint to the sheafification functor. This adjointness ensures  that
$F\mapsto\underline{H}^0(F)$ sends injective objects to injective objects.
Furthermore,  $F\mapsto \cH^p(U,F)$, $p\geq 0$ form  a universal $\delta$-functors from the
category of abelian presheaves on $U$ to the category of abelian groups,   see \cite{Tamme 1994}, Theorem 2.2.6.
In turn, we get a Grothendieck spectral sequence
\begin{equation}
\label{cech-to-sheaf}
\check{H}^p(U,\underline{H}^q(F))\Longrightarrow H^{p+q}(U,F),
\end{equation}
which we call \emph{\v{C}ech-to-sheaf-cohomology spectral sequence}.

\begin{theorem}
{\bf  (H.\ Cartan)}  
The following are equivalent:
\begin{enumerate}
\item $H^q(U,F)=0$ for all $U\in\shB$ and  $q\geq 1$.
\item $\check{H}^q(U,F)=0$ for all $U\in\shB$ and $q\geq 1$.
\end{enumerate}
If these equivalent conditions are satisfied, then $\check{H}^p(X,\underline{H}^q(F))=0$ for all $p\geq 0$ and $q\geq 1$,
and the edge maps $\check{H}^p(X,F)\ra H^p(X,F)$ are bijective for all $p\geq 0$.
\end{theorem}

\proof
Suppose (i) holds.
Then $\underline{H}^q(F)=0$ for all $q\geq 1$. This ensures that  $\check{H}^p(V,\underline{H}^q(F))=0$ for all $V\in\shC$.
In turn, the edge map
$\check{H}^p(V,F)\ra H^p(U,F)$ in the \v{C}ech-to-sheaf-cohomology spectral sequence (\ref{cech-to-sheaf}) is bijective.
For $V=U\in\shB$, this gives (ii), while for $V=X$ we get the amendment of the  statement.

Conversely, suppose that (ii) holds.
We inductively show that 
$H^p(U,F)=0$, or equivalently that the edge maps $\check{H}^p(U,F)\ra H^p(U,F)$ are bijective,   
for all $p\geq 1$ and all $U\in\shB$.
First note that $\check{H}^0(U,\underline{H}^r(F))=0$ for all $r\geq 1$,
because the sheafification of $\underline{H}^r(F)$ vanishes.
Thus the case $p=1$ in the induction   follows from the short exact sequence
$$
0\lra \check{H}^1(U,F)\lra H^1(U,F)\lra \check{H}^0(U,\underline{H}^1(F)).
$$
Now let $p\geq 1$ be arbitrary, and  suppose   that $H^p(U,F)=0$ for $1,\ldots, p-1$
and all $U\in\shB$. The argument in the previous paragraph gives $\check{H}^r(U,\underline{H}^s(F))=0$
for all $0<s<p$ and $r\geq 1$, and we already noted that $\check{H}^0(U,\underline{H}^p(F))=0$.
Consequently, the edge map $\check{H}^p(U,F)\ra H^p(U,F)$ must be bijective.
\qed

%===========================================================
\section{Appendix C:  Inductive dimension   in general topology}
\mylabel{InductiveDimension}

The proof of the crucial Proposition \ref{covering acyclic zariski} depends on   results
from dimension theory that are perhaps not so  well-known outside general topology.
Here we briefly review the relevant material. For a comprehensive treatment,
see for example Pears's monograph \cite{Pears 1975}, Chapter 4.

Let $X$ be a topological space. There are several ways to obtain   suitable notions
of dimension that are meaningful for compact spaces. 
One of them is the   \emph{large inductive dimension} $\Ind(X)\geq -1$, which goes  back to Brouwer \cite{Brouwer 1913}.
Here one first defines $\Ind(X)\leq n$ as a property of topological spaces  
by induction    as follows:
The induction starts with $n=-1$, where  $\Ind(X)\leq -1$ simply means that $X$ is empty.
Now suppose that  $n\geq 0$, and that the property is already defined for $n-1$.
Then $\Ind(X)\leq n$ means that for all closed subsets $A\subset X$
and all open neighborhoods $A\subset V$, there is a smaller open neighborhood
$A\subset U\subset V$ whose boundary 
$\Phi=\overline{U}\smallsetminus U$
has  $\Ind(\Phi)\leq n-1$. Now one defines $\Ind(X)=n$ as the least integer $n\geq -1$
for which $\Ind(X)\leq n$ holds, or $n=\infty$ if no such integer exists.

A variant is the \emph{small inductive dimension} $\ind(X)\geq -1$, which is defined in an analogous
way, but with points $a\in X$ instead of  closed subsets $A\subset X$.
If every point $a\in X$ is closed, an easy induction gives
$$
\ind(X)\leq \Ind(X).
$$
There are, however, compact spaces with strict inequality (\cite{Pears 1975}, Chapter 8, \S3).
Moreover, little useful can be said if the space   contains non-closed points.
For example, a finite linearly ordered space
$X=\left\{x_0,\ldots,x_n\right\}$  obviously has $\Ind(S)=0$ and  $\ind(S)=n$. Note that such spaces
arise as spectra  of valuation rings $R$ with Krull dimension $\dim(R)=n$.

Recall that a  space $X$ is  called \emph{at most zero-dimensional} if its topology admits
a basis consisting of open-and-closed subsets $U\subset X$. 
It is immediate that this is equivalent to $\ind(X)\leq 0$.
In fact, both 
small and large inductive dimension can be seen
as generalizations of this concept (\cite{Pears 1975}, Chapter 4, Proposition 1.1 and Corollary 2.2),
at least for compact spaces:

\begin{proposition}
\mylabel{CharacterizationZeroDimensional}
If $X$ is compact, then the following three conditions are equivalent: {\rm (i)}  $X$ is at most zero-dimensional;
{\rm (ii)} $\ind(X)\leq 0$;
{\rm (iii)} $\Ind(X)\leq 0$.
\end{proposition}

There are several \emph{sum theorems}, which express the large inductive dimension of a covering
$X=\bigcup A_\lambda$ in terms of the large inductive dimension of the subspaces $A_\lambda\subset X$.
The following sum theorem is very useful for us; its proof
is tricky, but only relies only   basic notations of set-theoretical topology
(\cite{Pears 1975}, Chapter 4, Proposition 4.13):

\begin{proposition}
\mylabel{SumTheorem}
Suppose that $X$ is compact, and $X=A\cup B$ is a   closed covering.
If the intersection $A\cap B$ is at most zero-dimensional, and the subspaces have  $\Ind(A),\Ind(B)\leq n$, 
then $\Ind(X)\leq n$.
\end{proposition}

From this we deduce the following fact, which, applied to $X_i=\Max(Y_i)$, enters in the proof
of Proposition \ref{covering acyclic zariski}.

\begin{proposition}
\mylabel{ZeroDimensionalCovering}
Suppose that $X$ is compact, and that $X=X_1\cup\ldots\cup X_n$ is a closed covering,
where the  $X_i$ are at most zero-dimensional. Then $X$ is at most zero-dimensional.
\end{proposition}

\proof
By induction it suffices to treat the case $n=2$.  
Since $X_i$ are at most zero-dimensional, the same holds for the subspace $X_1\cap X_2\subset X_i$.
The spaces $X_i$ are compact, because they are closed inside the compact space $X$.
According to Proposition \ref{CharacterizationZeroDimensional}, we have $\Ind(X_i)\leq 0$.
The sum theorem above yields $\Ind(X)\leq 0$, and Proposition \ref{CharacterizationZeroDimensional}  again tells us that
$X$ is at most zero-dimensional.
\qed

%===========================================================


\begin{thebibliography}{ccccc}

\bibitem{Aberbach 2005}
I.\ Aberbach:
The vanishing of $\Tor_1(R^+,k)$ implies that $R$ is regular. 
Proc.\ Amer.\ Math.\ Soc.\ 133 (2005),  27--29.

\bibitem{Artin 1970}
M.\ Artin:
Algebraization of formal moduli II: Existence of modifications.
Ann.\ Math.\  91 (1970), 88--135.

\bibitem{Artin 1971}
M.\ Artin:
On the joins of Hensel rings.
Advances in Math.\ 7 (1971),  282--296.

\bibitem{SGA 4a}
M.\ Artin, A.\ Grothendieck, J.-L. Verdier (eds.):
Th\'eorie des topos et cohomologie \'etale des sch\'emas (SGA 4) Tome 1.
Springer, Berlin, 1972.

\bibitem{Benoist 2013}
O.\ Benoist:
Quasi-projectivity of normal varieties. 
Int.\ Math.\ Res.\ Not.\ IMRN  17 (2013), 3878--3885. 

\bibitem{Bhatt 2014}
B.\ Bhatt: A flat map that is not a directed limit of finitely presented flat maps.
Note. http://www-personal.umich.edu/$\sim$bhattb/papers.html.

\bibitem{Bhatt; Scholze 2013}
B.\ Bhatt, P.\ Scholze:
The pro-\'etale topology for schemes.
Preprint (2013), arXiv:1309.1198.

\bibitem{Bhatt; de Jong 2014}
B.\ Bhatt, A.\ de Jong:
Lefschetz for local Picard groups.
Ann.\ Sci. \'Ec.\ Norm.\ Sup\'er.  47 (2014),  833--849.

\bibitem{TG 1-4}
N.\ Bourbaki:
Topologie g\'en\'erale. Chapitres 1--4. 
Hermann, Paris, 1971.

\bibitem{Brouwer 1913}
L.\ Brouwer:
\"Uber den nat\"urlichen Dimensionsbegriff.
J.\ Reine Angew.\ Math.\ 142 (1913), 146--152. 

\bibitem{Carral; Coste 1983}
M.\ Carral, M.\ Coste:
Normal spectral spaces and their dimensions.
J.\ Pure Appl.\ Algebra 30 (1983),  227--235. 

\bibitem{Conrad 2007}
B.\ Conrad:
Deligne's notes on Nagata compactifications.  
J.\ Ramanujan Math.\ Soc.\ 22 (2007),  205--257. 
Erratum in J.\ Ramanujan Math.\ Soc.\ 24 (2009), 427–428. 

\bibitem{Deligne 2010}
P.\ Deligne: 
Le th\'eor\`eme de plongement de Nagata.  
Kyoto J.\ Math.\ 50 (2010), 661--670. 

\bibitem{Endler 1972}
O.\ Endler:
Valuation theory. 
Springer, New York, 1972.

\bibitem{Flenner; O'Carroll; Vogel 1999}
H.\ Flenner, L.\  O'Carroll, W.\ Vogel:
Joins and intersections.  
Springer, Berlin, 1999.

\bibitem{Gabber 1994}
O.\ Gabber:
Affine analog of the proper base change theorem. 
Israel J.\ Math.\ 87 (1994),  325--335.

\bibitem{EGA I}
A.\ Grothendieck, J.\ Dieudonn\'e:
\'El\'ements de g\'eom\'etrie alg\'ebrique I:
Le langage des sch\'emas.
Springer, Berlin, 1970.

\bibitem{SGA 1}
A.\ Grothendieck:
Rev\^etements \'etales et groupe fondamental (SGA 1).
Springer, Berlin, 1971.

\bibitem{Gillman; Jerison 1960}
L.\ Gillman, M.\ Jerison:
Rings of continuous functions. 
Van Nostrand,  Princeton, N.J., 1960.

\bibitem{EGA II}
A.\ Grothendieck:
\'El\'ements de g\'eom\'etrie alg\'ebrique II:
\'Etude globale \'el\'ementaire de quelques classes de morphismes.
Publ.\ Math., Inst.\ Hautes \'Etud.\ Sci.\ 8 (1961).

\bibitem{EGA IVa}
A.\ Grothendieck:
\'El\'ements de g\'eom\'etrie alg\'ebrique IV: \'Etude locale des
sch\'emas et des morphismes de sch\'emas.
Publ.\ Math., Inst.\ Hautes \'Etud.\ Sci.\ 20 (1964).

\bibitem{EGA IVc}
A.\ Grothendieck:
\'El\'ements de g\'eom\'etrie alg\'ebrique IV: \'Etude locale des
sch\'emas et des morphismes de sch\'emas.
Publ.\ Math., Inst.\ Hautes \'Etud.\ Sci.\  28 (1966).

\bibitem{EGA IVd}
A.\ Grothendieck:
\'El\'ements de g\'eom\'etrie alg\'ebrique IV: \'Etude locale des
sch\'emas et des morphismes de sch\'emas.
Publ.\ Math., Inst.\ Hautes \'Etud.\ Sci.\   32 (1967).

\bibitem{SGA 2}
A.\ Grothendieck:
Cohomologie locale des faisceaux coh\'erents et th\'eor\`emes de Lefschetz locaux et globaux (SGA 2).
North-Holland Publishing Company, Amsterdam, 1968.

\bibitem{SGA 4b}
M.\ Artin, A.\ Grothendieck, J.-L. Verdier (eds.):
Th\'eorie des topos et cohomologie \'etale (SGA 4) Tome 2.
Springer, Berlin, 1973.

\bibitem{SGA 4c}
M.\ Artin, A.\ Grothendieck, J.-L. Verdier (eds.):
Th\'eorie des topos et cohomologie \'etale (SGA 4) Tome 3.
Springer, Berlin, 1973.

\bibitem{Borho; Weber 1971} 
W.\ Borho, V.\ Weber:
Zur Existenz total ganz abgeschlossener Ringerweiterungen. 
Math.\ Z.\ 121 (1971), 55--57. 

\bibitem{De Marco; Orsatti 1971}
G.\ De Marco, A.\ Orsatti:
Commutative rings in which every prime ideal is contained in a unique maximal ideal.
Proc.\ Amer.\ Math.\ Soc.\ 30 (1971) 459--466. 

\bibitem{Enochs 1968}
E.\ Enochs:
Totally integrally closed rings.
Proc.\ Amer.\ Math.\ Soc.\ 19 (1968), 701--706. 

\bibitem{Ferrand 2003}
D.\ Ferrand:
Conducteur, descente et pincement.
Bull.\ Soc.\ Math.\ France  131  (2003), 553--585.

\bibitem{Gillman and Jerison 1960}
L.\ Gillman, M.\ Jerison:
Rings of continuous functions. 
Van Nostrand,  Princeton, N.J., 1960.

\bibitem{Godement 1964}
R.\ Godement:
Topologie alg\'ebrique et  th\'eorie des faisceaux.
Hermann, Paris, 1964.

\bibitem{Gross 2012}
P.\ Gross:
The resolution property of algebraic surfaces.  
Compos.\ Math.\ 148 (2012),  209--226.

\bibitem{Hartshorne 1977}
R.\ Hartshorne:
Algebraic geometry.
Springer, Berlin,  1977.

\bibitem{Hochster 1970a}
M.\ Hochster:
Totally integrally closed rings and extremal spaces.
Pacific J.\ Math.\ 32 (1970), 767--779.

 \bibitem{Hochster 1970b}
M.\ Hochster:
Totally integrally closed rings and extremal spaces.
Pacific J.\ Math.\ 32 (1970), 767--779.

\bibitem{Hochster; Huneke 1992}
M.\ Hochster, C.\ Huneke:
Infinite integral extensions and big Cohen--Macaulay algebras.
Ann.\ of Math.\ 135 (1992), 53--89. 
 
\bibitem{Huneke; Lyubeznik 2007}
C.\ Huneke, G.\ Lyubeznik:
Absolute integral closure in positive characteristic.
Adv.\ Math.\ 210 (2007),  498--504. 

\bibitem{Huneke 2011}
C.\ Huneke:
Absolute integral closure. 
In:  A.\ Corso and C.\ Polini (eds.), Commutative algebra and its connections to geometry, pp.\ 119--135.
Amer. Math. Soc., Providence, RI, 2011. 

\bibitem{Horrocks 1971}
G.\ Horrocks:
Birationally ruled surfaces without embeddings in regular schemes.
Bull.\ London Math.\ Soc.\ 3 (1971), 57--60.

\bibitem{Johnstone 1977}
P.\ Johnstone:
Topos theory.
Academic Press, London, 1977.

\bibitem{Kleiman 1966}
S.\ Kleiman:
Toward a numerical theory of ampleness.
Ann.\ Math.\  84  (1966), 293--344.

\bibitem{Laumon; Moret-Bailly 2000}
G.\ Laumon, L.\ Moret-Bailly:
Champs alg\'ebriques.
Springer, Berlin, 2000.

\bibitem{Lazard 1967}
D.\ Lazard:
Disconnexit\'es des spectres d'anneaux et des pr\'esch\'emas
Bull.\ Soc.\ Math.\ France 95 (1967), 95--108.

\bibitem{Lenstra 1985}
H.\ Lenstra:
Galois theory for schemes.
Lecture notes, Universiteit Leiden, 1985.  
http://www.math.leidenuniv.nl/$\sim$hwl/PUBLICATIONS/pub.html.

\bibitem{Mazza; Voevodsky; Weibel 2006}
C.\ Mazza, V.\ Voevodsky, C.\  Weibel:
Lecture notes on motivic cohomology.
American Mathematical Society, Providence, RI, 2006. 

\bibitem{Milne 1980}
J.\ Milne: 
\'Etale cohomology.
Princeton University Press, Princeton, 1980.

\bibitem{Mumford 1961}
D.\ Mumford:
The topology of normal singularities of an algebraic surface and a criterion for simplicity.
Publ.\ Math., Inst.\ Hautes \'Etud.\ Sci.\ 9 (1961), 5--22.

\bibitem{Nicholson; Zhou 2005}
W.\ Nicholson, Y.\ Zhou:
Clean rings: a survey. 
In: J.\ Chen, N.\. Ding and H.\ Marubayashi (eds.), Advances in ring theory, pp.\ 181--198.
World Sci.\ Publ., Hackensack, NJ, 2005.

\bibitem{Nisnevich 1989}
Y.\ Nisnevich:
The completely decomposed topology on schemes and associated descent spectral sequences in algebraic K-theory. 
In: J.\ Jardine, V.\ Snaith (eds.), Algebraic K-theory: connections with geometry and topology, pp.\  241--342.
Kluwer, Dordrecht, 1989.

\bibitem{Pears 1975}
A.\ Pears:
Dimension theory of general spaces. 
Cambridge University Press, Cambridge, 1975. 

\bibitem{Perling; Schroeer 2014}
M.\ Perling, S.\ Schr\"oer:
Vector bundles on  proper toric 3-folds and certain other   schemes.
arXiv:1407.5443, to appear in Trans.\ Amer.\ Math.\ Soc.\

\bibitem{Rydh 2010}
D.\ Rydh:
Submersions and effective descent of \'etale morphisms. 
Bull.\ Soc.\ Math.\ France 138 (2010),  181--230. 

\bibitem{Rydh 2015}
D.\ Rydh:
Noetherian approximation of algebraic spaces and stacks.
J.\ Algebra 422 (2015), 105--147.

\bibitem{Schoutens 2004}
H.\ Schoutens:
On the vanishing of Tor of the absolute integral closure. 
J.\ Algebra 275 (2004),  567--574. 

\bibitem{Schroeer 1999}
S.\ Schr\"oer:
On non-projective normal surfaces.
Manuscr.\ Math.\ 100 (1999), 317--321.

\bibitem{Schroeer 2014}
S.\ Schr\"oer:
Points in the fppf topology.
arXiv:1407.5446.
To appear in  Ann.\ Sc.\ Norm.\ Super.\ Pisa Cl.\ Sci.

\bibitem{Schwartz 2011}
N.\ Schwartz:
Graph components of prime spectra.  
J.\ Algebra 337 (2011), 13--49.

\bibitem{Sergio 1985}
I.\ Sergio:
Rings whose spectra do not have polygons. 
Matematiche (Catania) 37 (1982), 62--71 (1985).

\bibitem{stacks project}
Stacks project contributors: Stacks project.
http://stacks.math.columbia.edu/browse,
version 2676bbe, compiled on Jun 17, 2014.

\bibitem{Tamme 1994}
G.\ Tamme:
Introduction to \'etale cohomology. 
Springer, Berlin, 1994.

\bibitem{Thomason; Trobaugh 1990}
R.\ Thomason, T.\ Trobaugh:
Higher algebraic $K$-theory of schemes and of derived categories. 
In: P.\ Cartier et al.\ (eds.),
The Grothendieck Festschrift III, 247--435.
Birkh\"auser, Boston, 1990.

\bibitem{Vakil and Wickelgren 2011}
R.\ Vakil, K.\ Wickelgren:
Universal covering spaces and fundamental groups in algebraic geometry as schemes.  
J.\ Th\'eor.\ Nombres Bordeaux 23 (2011), 489--526. 

\bibitem{Voevodsky; Suslin; Friedlander 2000}
V.\ Voevodsky, A.\ Suslin, E.\ Friedlander:
Cycles, transfers, and motivic homology theories.
Princeton University Press, Princeton, NJ, 2000.

\end{thebibliography}
\end{document}